\documentclass[journal]{IEEEtran}
\usepackage{graphicx}
\usepackage{amsmath}
\usepackage{amssymb}
\usepackage{CJK}
\usepackage{color}
\usepackage{dsfont}
\usepackage{subfigure}
\usepackage{booktabs}
\usepackage{microtype}
\usepackage{mathrsfs}
\usepackage{bm}
\usepackage{cite}
\usepackage{url}
\usepackage[colorlinks=true,allcolors=blue]{hyperref}

\ifCLASSINFOpdf
\else
\fi

\begin{document}
\title{Data-Based Distributionally Robust \\ Stochastic Optimal Power Flow -- Part II: \\ Case Studies}

\author{Yi Guo,~\IEEEmembership{Student Member,~IEEE,}
        Kyri Baker,~\IEEEmembership{Member,~IEEE,}
        Emiliano Dall'Anese,~\IEEEmembership{Member,~IEEE,}\\
        Zechun Hu,~\IEEEmembership{Senior Member,~IEEE,}
        Tyler H. Summers,~\IEEEmembership{Member,~IEEE}% <-this % stops a space
\thanks{This material is based on work supported by the National Science Foundation under grant CNS-1566127. \emph{(Corresponding author: Tyler H. Summers)}}
\thanks{Y. Guo and T.H. Summers are with the Department
of Mechanical Engineering, The University of Texas at Dallas, Richardson,
TX, USA, email: \{yi.guo2,tyler.summers\}@utdallas.edu.}
\thanks{K. Baker is with the Department of Civil, Environmental, and Architectural Engineering, The University of Colorado, Boulder, CO, USA, email:kyri.baker@colorado.edu.}
\thanks{E. Dall'Anese is with the Department of Electrical, Computer, and Energy Engineering, University of Colorado Boulder, Boulder, CO, USA, email: emiliano.dallanese@colorado.edu.}
\thanks{Z. Hu is with the Department of Electrical Engineering, Tsinghua University, Beijing, China, email: zechhu@tsinghua.edu.cn.}
}

\maketitle
\begin{abstract}
 This is the second part of a two-part paper on data-based distributionally robust stochastic optimal power flow (OPF). The general problem formulation and methodology have been presented in Part I \cite{Guo2018DataDriven1}. Here, we present extensive numerical experiments in both distribution and transmission networks to illustrate the effectiveness and flexibility of the proposed methodology for balancing efficiency, constraint violation risk, and out-of-sample performance. On the distribution side, the method mitigates overvoltages due to high photovoltaic penetration using local energy storage devices. On the transmission side, the method reduces N-1 security line flow constraint risks due to high wind penetration using reserve policies for controllable generators. In both cases, the data-based distributionally robust model predictive control (MPC) algorithm explicitly utilizes forecast error training datasets, which can be updated online. The numerical results illustrate inherent tradeoffs between the operational costs, risks of constraints violations, and out-of-sample performance, offering systematic techniques for system operators to balance these objectives.
\end{abstract}


% \fi
\IEEEpeerreviewmaketitle

\section{Introduction}


% background, motivation

% context
\IEEEPARstart{A}s penetration levels of renewable energy sources (RESs) continue increasing to substantial fractions of total supplied power and energy, system operators will require more sophisticated methods for balancing inherent tradeoffs between nominal performance and operational risks. This is relevant for both transmission and distribution networks, which have seen rapid recent increases in photovoltaic (PV) and wind energy sources. Here, we present extensive numerical experiments in both distribution and transmission networks to illustrate the effectiveness and flexibility of our data-based distributionally robust stochastic optimal power flow (OPF) methodology for balancing efficiency, constraint violation risk, and out-of-sample performance.

% brief lit review
A variety of stochastic OPF methods have been recently proposed to explicitly incorporate forecast errors of network uncertainties for modeling risk. These methods can be categorized according to assumptions made about forecast error distributions and metrics for quantifying risk. Many make specific assumptions about forecast error distributions and use chance constraints, which encode value at risk (VaR), conditional value at risk (CVaR), or distributional robustness with specific ambiguity sets \cite{bienstock,zhang1,perninge,roald1,tyler,LiMathieu,zhang2,hedayati2018reserve,bazrafshan2014voltage}. Others handle uncertainty using scenarios approaches, sample average approximation, or robust optimization \cite{yong,capitanescu,conejo,Vrakopoulou1,Vrakopoulou2,warrington,jabr,lubin ,roald2016integrated,DallAnese2017chance,hug,hug2008coordinated,vrakopoulou3,li2017chance,li2017distribued,li2017distribution,mohagheghi2017framework,mohagheghi2018real}. However, in practice all decisions about how to model forecast error distributions must be based on finite historical training datasets. None of the existing methods account for sampling errors inherent in such datasets in high-dimensional spaces.

% contributions
Our proposed data-based distributionally robust stochastic OPF methodology, developed in Part I \cite{Guo2018DataDriven1},  explicitly incorporates robustness to sampling errors by considering sets of distributions centered around a finite training dataset. A model predictive control (MPC) approach utilizes Wasserstein-based distributionally robust optimization subproblems to obtain superior out-of-sample performance. Here in Part II we illustrate the effectiveness and flexibility of the methodology in both distribution and transmission networks. On the distribution side, the method mitigates overvoltages due to high photovoltaic penetration using local energy storage devices. On the transmission side, the method reduces N-1 security line flow constraint risks due to high wind penetration using reserve policies for controllable generators. In both cases, the data-based distributionally robust (MPC) algorithm explicitly utilizes forecast error training datasets, which can be updated online. The numerical results illustrate inherent tradeoffs between the operational costs, risks of constraints violations, and out-of-sample performance, offering systematic techniques for system operators to balance these objectives.

The rest of the paper is organized as follows. Section II illustrates the data-based distributionally robust stochastic AC OPF for mitigating overvoltages in a modified IEEE 37-node  distribution feeder with high PV penetration using local energy storage devices. Section III illustrates the data-based distributionally robust stochastic DC OPF for reducing N-1 security line flow constraint risks due to high wind penetration using reserve policies for controllable generators. Section IV concludes.
\\

\textbf{Notation}: The inner product of two vectors $a,b \in \mathbf{R}^m$  is denoted by $\langle a, b \rangle := a^\intercal b$. The $N_s$-fold product of distribution $\mathds{P}$ on a set $\Xi$ is denoted by $\mathds{P}^{N_s}$, which represents a distribution on the Cartesian product space $\Xi^{N_s} = \Xi \times \ldots \times \Xi$. We use $N_s$ to represent the number of samples inside the training dataset $\hat{\Xi}$. Superscript `` $\hat{\cdot}$ " is reserved for the objects that depend on a training dataset $\hat{\Xi}^{N_s}$. We use $(\cdot)^\intercal$ to denote  vector or matrix transpose.
%and $A^{\mathcal{H}_t}$ to denote matrix $A$ to the power $\mathcal{H}_t$, where $\mathcal{H}_t$ is a finite time horizon. Let $\mathcal{T}_t$ indicate a time interval $[t, t+\mathcal{H}_t]$, which is also written in $[\underline{t},\bar{t}]$ for mathematical convenience.
The operators $\Re\{\cdot\}$ and $\Im\{\cdot\}$ return the real and imaginary part of a complex number, respectively. The operator $[\,\cdot\,]_{[a,b]}$ selects the $a$-th to $b$-th elements of a vector or rows of a matrix.

% data-based stochastic opf in distribution networks.
\section{Overvoltage mitigation in distribution networks}
In this section, we apply the data-based distributionally robust stochastic OPF methodology to mitigate overvoltages  in distribution networks by controlling set points in RESs and energy storage devices. We provide further modeling details of the loads, inverter-based RESs, and energy storage devices. The set points of controllable devices are repeatedly optimized over a finite planning horizon within a MPC feedback scheme. The risk conservativeness of the voltage magnitude constraints and the out-of-sample performance robustness to sampling errors are explicitly adjustable by two scalar parameters.

\subsection{System model}
\subsubsection{Loads} We use $P_{l,n}^t$ and $Q_{l,n}^t$ to denote the active and reactive power demands at bus $n \in \mathcal{N}$. We also define two vectors $\mathbf{p}_l^t :=[P_{l,1}^t, \ldots ,P_{l,N}^t]^\intercal$ and $\mathbf{q}_l^t :=[Q_{l,1}^t, \ldots, Q_{l,N}^t]^\intercal$. If no load is connected to bus $n \in \mathcal{N}$, then $P_{l,n}^t = 0$ and $Q_{l,n}^t = 0$. Load uncertainties are modeled based on historical data of forecast errors. The active and reactive loads are given by $\mathbf{p}_l^t = \bar{\mathbf{p}}_l^t(u_t) + \tilde{\mathbf{p}}_l^t(\xi_t)$, $\mathbf{q}_l^t = \bar{\mathbf{q}}_l^t(u_t) + \tilde{\mathbf{q}}_l^t(\xi_t)$, where $\mathbf{\bar{p}}_l^t(u_t) \in \mathbf{R}^N$ and $\mathbf{\bar{q}}_l^t(u_t) \in \mathbf{R}^N$ are forecasted nominal loads, which can depend on control decisions (e.g., load curtailment control). The nodal injection errors $\tilde{\mathbf{p}}_l^t(\xi_t) \in \mathbf{R}^N$ and $\tilde{\mathbf{q}}_l^t(\xi_t) \in \mathbf{R}^N$ depend on the aggregate forecast error vector $\xi_t$.

\subsubsection{Renewable energy model} Let $P_{\textrm{av},n}^t$ be the maximum availability renewable energy generation at bus $n \in \mathcal{N}_R \subseteq \mathcal{N} $, where the set $\mathcal{N}_R$ denotes all buses with RESs. With high RES penetration, overvoltages can cause power quality and reliability issues. By intelligently operating set points of RES and energy storage, operators can optimally trade off risk of constraint violation and economic efficiency (e.g., purchase of electricity from the main grid, active power curtailment costs, and reactive compensation costs). The active power injections of RESs are controlled by adjusting an active power curtailment factor $\alpha_n^t \in [0,1]$. Reactive power set points of RESs can also be adjusted within a limit $\bar{S}_n$ on apparent power as follows
\begin{equation} \nonumber
\sqrt{((1-\alpha_n^t)P_{\textrm{av},n}^t)^2 + (Q_n^t)^2} \leq \bar{S}_n, n \in \mathcal{N}_R.
\end{equation}
We define aggregate vectors: $\bm{\alpha}^t :=[\alpha_1^t, \ldots ,\alpha_N^t]^\intercal$ and $\mathbf{p}_{\textrm{av}}^t:= [P_{\textrm{av},1}^t,\ldots, P_{\textrm{av},N}^t]^\intercal$ and $\mathbf{q}_{c}^t:=[Q_{1}^t,\ldots, Q_{N}^t]^\intercal$.

If bus $n \in \mathcal{N} \backslash \mathcal{N}_R$ has no RES, by convention we set $\alpha_n^t = 0, P_{\textrm{av},n}^t = 0$ and $Q_n^t = 0$. The curtailment factor and reactive power compensation $\{\alpha_{n}^t$, $Q_{n}^t\}$ together set the inverter operating point and are also subject to a power factor constraint
\begin{equation*}
 |Q_n^t| \leq \textrm{tan}(\theta_n)[(1-\alpha_n^t)P_{\textrm{av},n}^t], n\in \mathcal{N}_R,
\end{equation*}
where $\textrm{cos}(\theta_n) \in (0,1]$ is the power factor limit for RESs. The power factor constraint is convex, and can be discarded in settings where the inverters are not required to  operate at a minimum power factor level. The premise here is that RESs can assist in the regulation of voltages by promptly adjusting the reactive power and curtailing active power as needed; RESs can provide faster voltage regulation capabilities compared to traditional power factor correction devices (i.e., capacitor banks). The proposed data-based distributionally robust OPF will consider adjustments of both active and reactive powers to aid voltage regulation, which in principle can be done in both transmission and distribution networks.

% Energy Storage Model
\subsubsection{Energy storage model} The state-of-charge (SOC) of the energy storage device located at bus $n\in\mathcal{N}_B \subseteq \mathcal{N}$ in kWh is represented as $B_n^t$. The dynamics of these devices are
\begin{equation}\label{batterdynamic}
B_n^{t+1} = B_n^t + \eta_{B,n} P_{B,n}^t \Delta, n \in \mathcal{N}_B,
\end{equation}
where $\Delta$ is the duration of the time interval $(t,t+1]$, and $P_{B,n}^t$ is the charging/discharging power of the storage device in kW. We assume the battery state is either charging ($P_{B,n}^t \geq 0$) or discharging ($P_{B,n}^t \leq 0$) during each time interval $(t, t+1]$. For simplicity, we suppose the round-trip efficiency of the storage device $\eta_{B,n}=1$ to avoid the nonconvexity when introducing binary variables. Additionally, two common operational constraints of energy storage devices are
\begin{equation}\nonumber
B_n^{\textrm{min}} \leq B_{n}^t \leq B_n^{\textrm{max}},~~~~
P_{B,n}^{\textrm{min}} \leq P_{B,n}^t \leq P_{B,n}^{\textrm{max}},
\end{equation}
where $B_n^{\textrm{min}}$, $B_n^{\textrm{max}}$ are the rated lower and upper SOC levels, and $P_{B,n}^{\textrm{min}}$, $P_{B.n}^{\textrm{max}}$ are the minimum and maximum charging/discharging limits. Other constraints can be added for electric vehicles (EVs), for example, a prescribed SOC $B_{n}^t = B_n^{\textrm{max}}$ at a particular time. If no energy device is connected to a certain bus, the charging/discharging power and SOC are fixed to zero: $P_{B,n}^t = 0$, $B_n^t =0$, for all $n \in \mathcal{N}\backslash\mathcal{N}_B$. We define the aggregate vectors $\mathbf{p}_B^t := [P_{B,1}^t,\ldots,P_{B,N}^t]^\intercal$, and $\mathbf{b}^t := [B_1^t,\ldots,B_N^t]^\intercal$.

\subsection{Data-based stochastic OPF implementation}
We now use the methodologies presented in Part I Section IV, and the models of loads, RESs and energy storage devices to develop a data-based stochastic AC OPF for solving a voltage regulation problem. This stochastic OPF aims to balance the  operational costs the total CVaR values of the voltage magnitude constraints. We consider an operational cost that captures electricity purchased by customers, excessive solar energy fed back to the utility, reactive power compensation costs and penalties for active power curtailment
\begin{equation}\nonumber
\begin{split}
	& J_\textrm{Cost}^t (\bm{\alpha}^t, \mathbf{q}_c^t, \mathbf{p}_B^t, \xi_t) =\\
	& = \sum_{n\in\mathcal{N}}a_{1,n}^t \left[ P_{l,n}^t + P_{B,n}^t - (1-\alpha_n^t)P_{\textrm{av},n}^t \right]_+\\
	& + \sum_{n\in\mathcal{N}}a_{2,n}^t \left[(1-\alpha_n^t)P_{\textrm{av},n}^t - P_{l,n}^t - P_{B,n}^t\right]_+\\
	& + \sum_{n\in\mathcal{N}}a_{3,n}^t |Q_n^t| + \sum_{n\in\mathcal{N}}a_{4,n}^t\alpha_n^t P_{\textrm{av},n}^t .
\end{split}
\end{equation}
We collect all decision variables into $\mathbf{y}_t = \{\alpha^t, \mathbf{q}_c^t, \mathbf{p}_B^t, \mathbf{b}^t\}$, and all RES and load forecast errors into the random vector $\xi_t$. Now the MPC subproblems take the following form
\\
\textbf{Data-based distributionally robust stochastic OPF}
\begin{subequations} \label{DBACOPF}
\begin{align}
	&\hspace{-2mm}\inf_{\begin{subarray}{c} \mathbf{y}_\tau, \kappa_o^\tau, \\ \varpi_{1,n}, \varpi_{2,n} \end{subarray}} \sum_{\tau = t}^{t + \mathcal{H}_t} \bigg\{\mathds{E} [\hat{J}^\tau_{\textmd{Cost}}] +  \rho \hspace{-2 mm} \sup_{\mathds{Q}_\tau \in \hat{\mathcal{P}}^{N_s}_\tau}  \sum_{o=1}^{N_\ell} \mathds{E}^{\mathds{Q}_\tau}[\bar{\mathcal{Q}}_o^\tau]  \bigg\}, \nonumber\\
	&\hspace{-2mm}=\hspace{-1mm} \inf_{\begin{subarray}{c}  \mathbf{y}_\tau, \kappa_o^\tau ,\\ \varpi_{1,n}, \varpi_{2,n}\\
    \lambda_o^\tau, s_{io}^\tau, \varsigma_{iko}^\tau \end{subarray}} \sum_{\tau=t}^{t + \mathcal{H}_t} \bigg\{ \mathds{E}[\hat{J}^\tau_{\textrm{Cost}}] + \sum_{o=1}^{N_\ell} \bigg(\lambda_o\varepsilon_\tau + \frac{1}{N_s}\sum_{i=1}^{N_s} s_{io}^\tau\bigg) \bigg\},\\
	&\hspace{-2mm}\textrm{subject to}\nonumber\\
	&\hspace{-2mm}\rho (\bar{\mathbf{b}}_{ok}(\kappa_o^\tau) + \langle \bar{\mathbf{a}}_{ok}(\mathbf{y}_\tau), \hat{\xi}_\tau^{i}\rangle + \langle\varsigma_{iko}, \mathbf{d}-H\hat{\xi}_\tau^{i}\rangle) \le s_{io}^\tau, \label{DDRO1} \\ 
	&\hspace{-2mm}\| H^\intercal \varsigma_{iko}-\rho \bar{\mathbf{a}}_{ok}(\mathbf{y}_\tau)\|_\infty \le \lambda_o^\tau,\label{DDRO2}\\ 
	&\hspace{-2mm}\varsigma_{iko} \geq 0,\label{DDRO3}\\ 
	&\hspace{-2mm}\frac{1}{N_s}\sum_{i=1}^{N_s}\bigg[[(1-\alpha^\tau_n) \hat{P}_{\textrm{av},n}^{\tau,i}]^2 + (Q_n^\tau)^2 - \bar{S}_n^2 + \varpi_{1,n}^\tau \bigg]_+ \leq \varpi_{1,n}^\tau \beta, \label{inverterspowerconstraints} \\
	&\hspace{-2mm}\frac{1}{N_s}\sum_{i=1}^{N_s}\bigg[\textrm{tan}(\theta_n)[(1-\alpha_n^t)\hat{P}_{\textrm{av},n}^{\tau,i}] - |Q_n^\tau| + \varpi_{2,n}^\tau \bigg]_+ \leq \varpi_{2,n}^\tau \beta, \label{powerfactor}\\
	&\hspace{-2mm}B_n^{\textrm{min}} \leq B_{n}^\tau \leq B_n^{\textrm{max}},\\
	&\hspace{-2mm}P_{B,n}^{\textrm{min}} \leq P_{B,n}^\tau \leq P_{B,n}^{\textrm{max}},\\
	&\hspace{-2mm}B_n^{\tau+1} = B_n^\tau + \eta_{B,n} P_{B,n}^\tau \Delta, \\
	&\hspace{-2mm}0 \leq \alpha_n^\tau \leq 1,\\
	&\hspace{-2mm}\forall i\le N_s, ~\forall o \le N_\ell, ~n \in \mathcal{N}_R, ~k=1,2, ~\tau = t,\ldots, t+\mathcal{H}_t,\nonumber
\end{align}
\end{subequations}
where $\varpi_{1,n}^\tau$, $\varpi_{2,n}^\tau$, and $\kappa_o^\tau$ are CVaR auxiliary valuables, and $\lambda_o^\tau$, $s_{io}^\tau$, $\varsigma_{iko}^\tau$ are auxiliary variables associated with the distributionally robust Wasserstein ball reformulation. For simplicity, the power factor constraints and apparent power limitation constraints are not treated as distributionally robust constraints, and instead are handled using direct sample average approximation.
% Simulation Results Start from Here.

\textbf{Remark 2.1} (battery efficiency). To maintain convexity of the underlying problem formulation and therefore facilitate the development of computationally affordable solution methods, we utilized an approximate model for the battery dynamics with no charging and discharging efficiency losses  \eqref{DBACOPF}. At the expense of significantly increasing the problem complexity, charging and discharging efficiencies can be accommodated as \cite{rahbar2015real}
\begin{equation*}
B_n^{\tau +1} = B_n^\tau + \eta_c P_{B_c,n}^\tau\Delta - \frac{1}{\eta_{d}}P_{B_d,n}^\tau\Delta,
\end{equation*} 
where $\eta_c, \eta_d \in (0,1]$ denote the charging and the discharging efficiencies, respectively; $P_{B_c,n}^\tau \geq 0$ represents the charging rate and $P_{B_d,n}^\tau \geq 0$ the discharging rate at time $\tau$. Additional constraints, however, are needed to ensure that the solution avoids meaningless solutions where a battery is required to charge and discharge simultaneously; in particular, one can: a) add a constraint $P_{B_c,n}^\tau P_{B_d,n}^\tau = 0$ \cite{rahbar2015real}; or. b) introduce binary variables to indicate the charging status (e.g., charging/discharging) of the batteries \cite{jabr}. Either way, given the non-convexity of the resultant problem, possibly sub-optimal solutions can be achieved~\eqref{DBACOPF}. In addition, exact relaxation methods under appropriate assumptions offer an alternative way to maintain convexity of the charging problem; see \cite{li2016sufficient,ding2015value}. Extending the proposed technical approach to a setting with binary variables or exact relaxation methods will be pursued as a future research effort. 

\textbf{Remark 2.2} (battery life). The degradation of energy storage systems may depend on the depth of discharge and  the number of charging/discharging cycles \cite{divya2009battery}. The battery aging process is usually described by partial differential equations \cite{ramadesigan2012modeling}; this is a practical model for industrial applications,  but it introduces  significant computational challenges in optimization tasks \cite{xu2018optimal}. Additional optimization variables as well as penalty functions could be included to limit the number of cycles per day and ensure a minimum state of charge~\cite{baker2017energy,fortenbacher2017modeling,ortega2014optimal,ying2016stochastic,koller2013defining,shi2018using,zhang2017evaluation}. Pertinent reformulations to account for battery degradation will be pursued in future research activities. 

\textbf{Remark 2.3} (voltage at slack bus). Similar to the majority of the works in the literature, the  voltage at the slack bus (i.e., substation) is considered as an input of the problem (and, therefore, it is not controllable).   However, it is worth noting that discrete variables modeling  changing the tap position of the transformer   can be  incorporated in \eqref{DBACOPF}; see e.g., \cite{mohagheghi2018real,mohagheghi2017framework}. Branch and bound techniques can then be utilized to solve the problem. 

\subsection{Numerical results}
\begin{figure}[htbp!]
\centering
\includegraphics[width=2.8in]{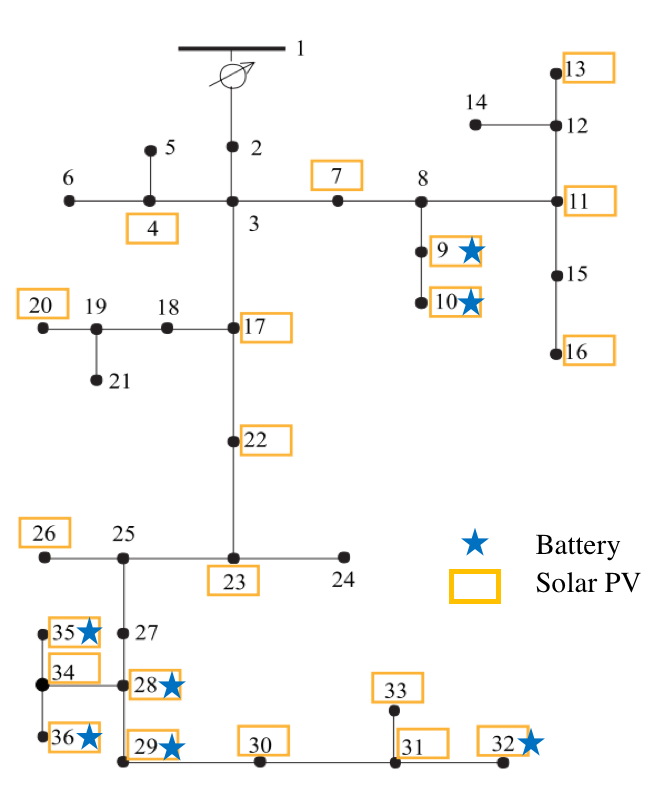}
\caption{IEEE 37-node test feeder with renewable energy resources and storage devices.}
\label{fig:distributionfeeder}
\end{figure}

\begin{figure}[t!]
\centering
\includegraphics[width=3.5in]{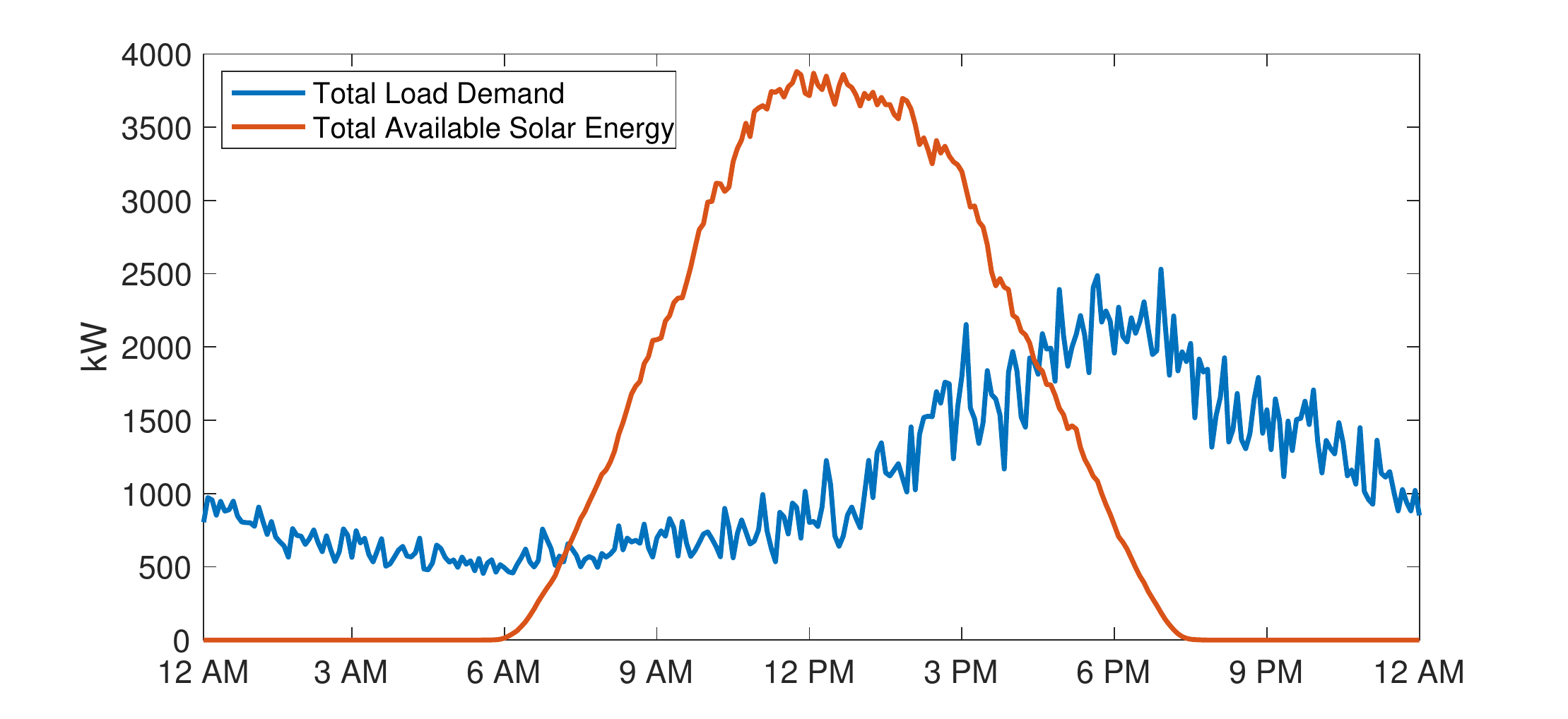}
\caption{Total available solar energy and load demand.}
\label{fidg:loadsolardata}
\end{figure}

We use a modified IEEE-37 node test feeder to demonstrate our proposed data-based stochastic AC OPF method. As shown in Fig. \ref{fig:distributionfeeder}, the modified network is a single-phase equivalent and the load data is derived from real measurements from feeders in Anatolia, CA during the week of August 2012 \cite{bank2013development}. We place 21 photovoltaic (PV) systems in the network. Their locations are marked by yellow boxes in Fig.  \ref{fig:distributionfeeder}, and their capacities are summarized in Table \ref{table:capacityBatteryRES}. Based on irradiation data from \cite{bank2013development,solardata}, we utilized a greedy gradient boosting method \cite{friedman2001greedy} to make multi-step ahead predictions of solar injections, and then computed a set of forecast errors from the dataset. In general forecast errors increase with the prediction horizon.
Other parameters of the network, such as line impedances and shunt admittances, are taken from \cite{zimmerman2011matpower}. The total nominal available solar power $\sum_{n}P_{\textrm{av},n}^t$ and aggregate load demand over 24 hours is also shown in Fig. \ref{fidg:loadsolardata}.

% a table for the solar energy data
\begin{table}[htbp!]
\caption{Capacities of inverted-based solar energy generations and energy storage devices}
\centering
\begin{tabular}{|c|c|c|c|c|c|}
\hline
Node & $S_n$ [kVA] & Node & $S_n$ [kVA] & Node & $B_n^{\rm{max}}$\\
\hline
4 & 150 & 7 & 300 & 9 & 100\\
\hline
9 & 300 & 10 & 600 & 10 & 100\\
\hline
11 & 660 & 13 & 360 & 28 & 50\\
\hline
16 & 600 & 17 & 360 & 29 & 250\\
\hline
20 & 450 & 22 & 150 & 32 & 250\\
\hline
23 & 750 & 26 & 300 & 35 & 120\\
\hline
28 & 750 & 29 & 300 & 36 & 200\\
\hline
30 & 360 & 31 & 600\\
\cline{1-4}
32 & 330 & 33 & 750\\
\cline{1-4}
34 & 450 & 35 & 450\\
\cline{1-4}
36 & 450\\
\cline{1-2}
\end{tabular}
\label{table:capacityBatteryRES}
\end{table}

The energy storage systems are placed with PV systems at certain nodes, as shown in Fig. \ref{fig:distributionfeeder}. Their locations and capacities are listed in Table \ref{table:capacityBatteryRES}. We select the capacities in the range of typical commercial storage systems, or aggregate of 10-12 residential-usage batteries (e.g., electric vehicles), which are connected to the same step-down transformer. The lower limit of SOC, $B_{n}^{\textrm{min}}$, is set to be zero for all batteries. The charging/discharging rate $P_{B,n}^t$ is also limited by 10\% of their respective energy capacity $B_{n}^{\textrm{max}}$. Voltage limits $V^{\rm{max}}$ and $V^{\rm{min}}$ are 1.05 p.u. and 0.95 p.u., respectively. The cost function parameters are $a_{1,n}^t = 10$, $a_{2,n}^t = 3$, and $a_{3,n}^t =3$ and $a_{4,n}^t = 6$. The decision making time period is 5 minutes.

Due to high PV penetration, overvoltage conditions can emerge during solar peak irradiation. Given the real data available for the numerical tests, the numerical tests are focused on alleviating over-voltage conditions via the proposed distributionally-robust tools. 
Other constraints are approximated via sample average methods~\cite{linderoth2006empirical,tyler,DallAnese2017chance}; however, in general it is straightforward to formulate other constraints as distributionally robust. The power factor $\textrm{PF}$ limit is 0.9 in \eqref{powerfactor}. The risk level parameter $\eta$ is set to 0.01 for quantifying 1\% violation probability of constraints \eqref{inverterspowerconstraints}-\eqref{powerfactor}. To emphasize the effect of sampling errors, the number of forecast error samples $N_s$ included in the training dataset $\hat{\Xi}_t^{N_s}$ is limited to 30. The forecast errors are not assumed bounded, so the parameters of the support polytope $\Xi_t := \{\xi_t \in \mathbf{R}^{N_\xi} : H\xi_t \leq \mathbf{d} \}$  are set to zero in \eqref{DDRO1}-\eqref{DDRO2}. We solved \eqref{DBACOPF} using the MOSEK  solver \cite{mosek} via the MATLAB interface CVX \cite{cvx} on a laptop with 16 GB of memory and 2.8 GHz Intel Core i7. Solving each time step during solar peak hours with distibutionally robust constraints takes 4.84 seconds. Note that our implementation is not optimized for speed and in principle could easily be sped up and scaled to larger problems since the problem is ultimately convex quadratic.
\color{black}

\begin{figure*}[htbt!]
\centering
\subfigure{\label{fig:PVcurtail1} %% label for first subfigure
\includegraphics[width=1.7in]{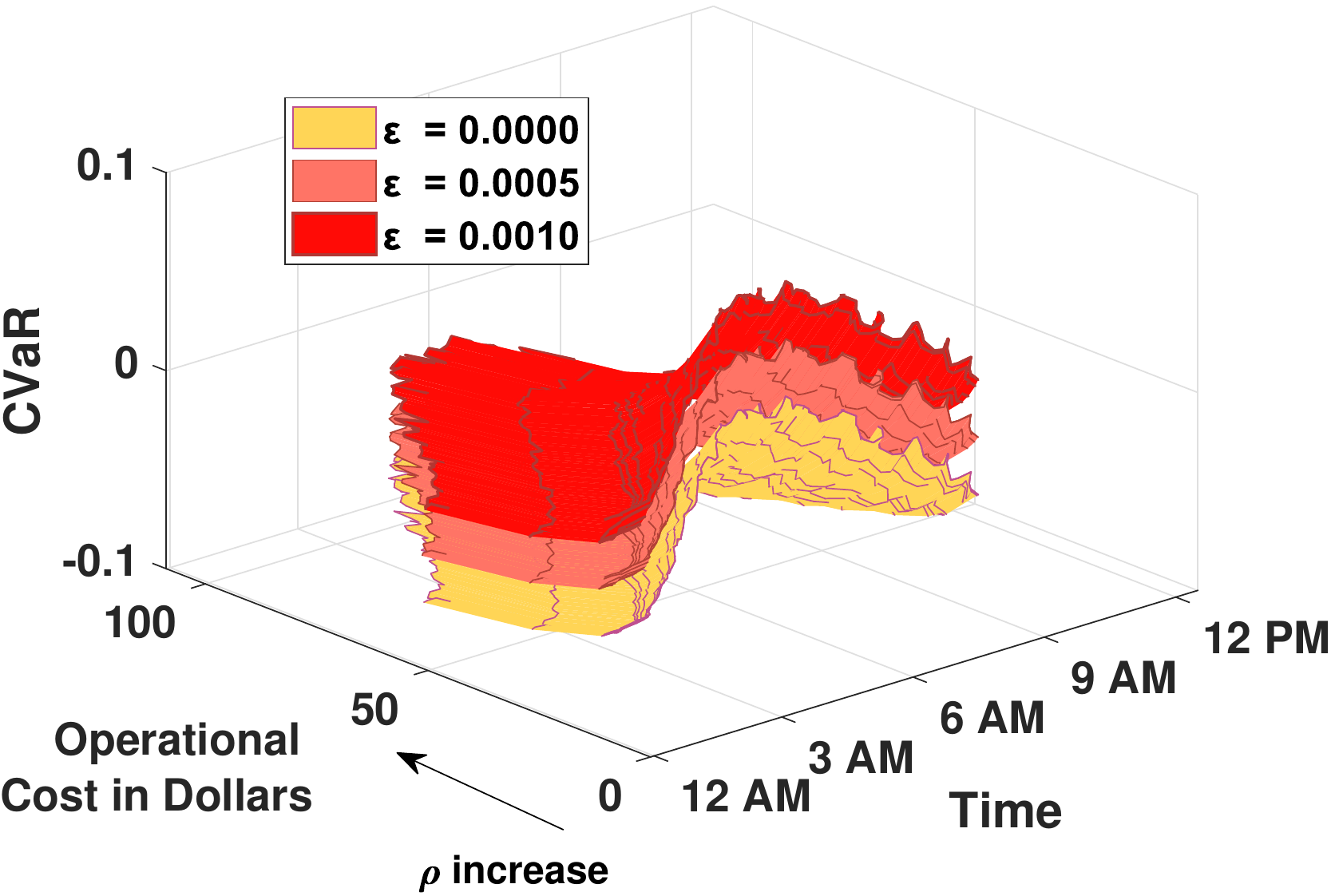}}
\subfigure{\label{fig:PVcurtail1} %% label for first subfigure
\includegraphics[width=1.7in]{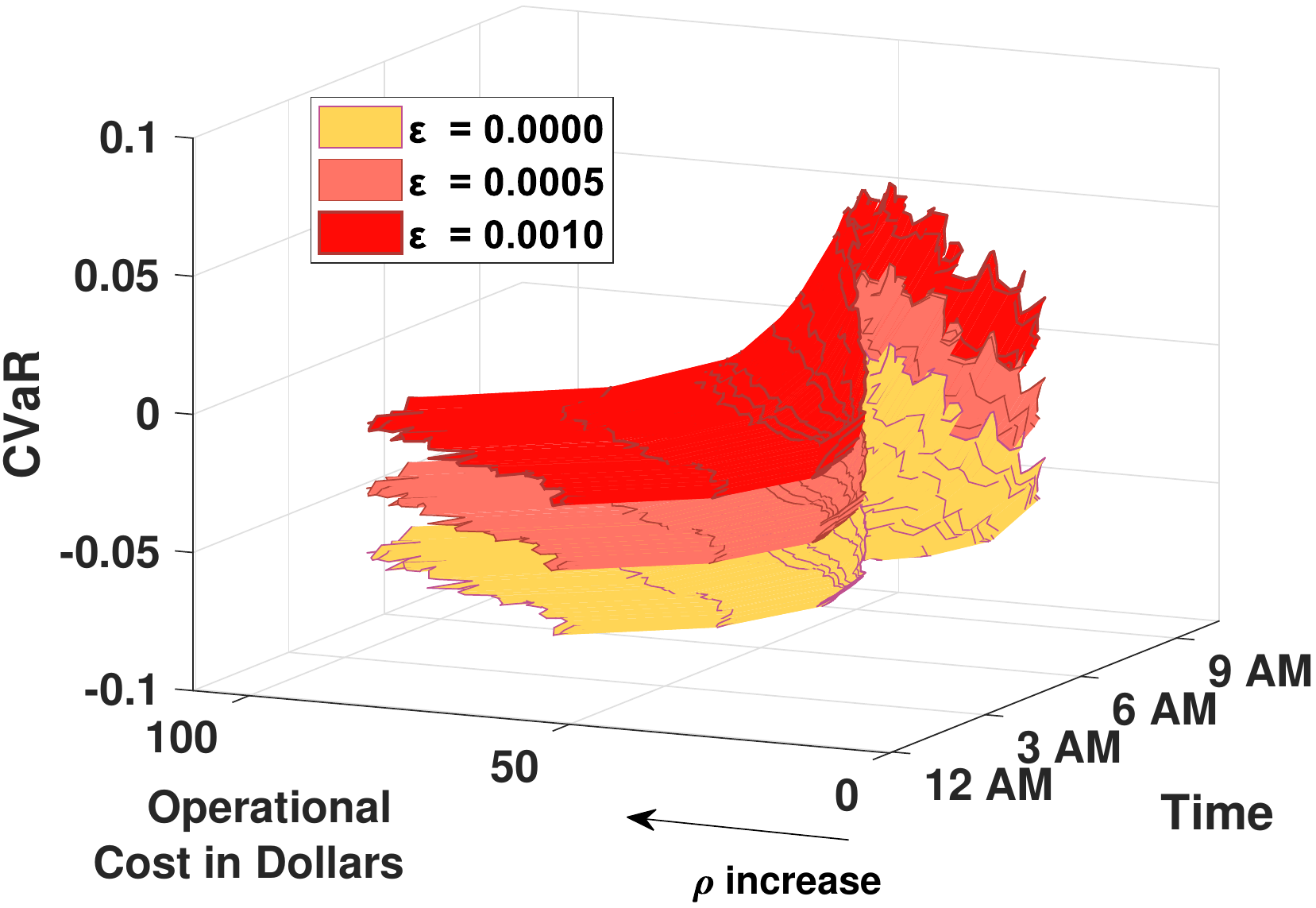}}
\subfigure{\label{fig:PVcurtail1} %% label for first subfigure
\includegraphics[width=1.7in]{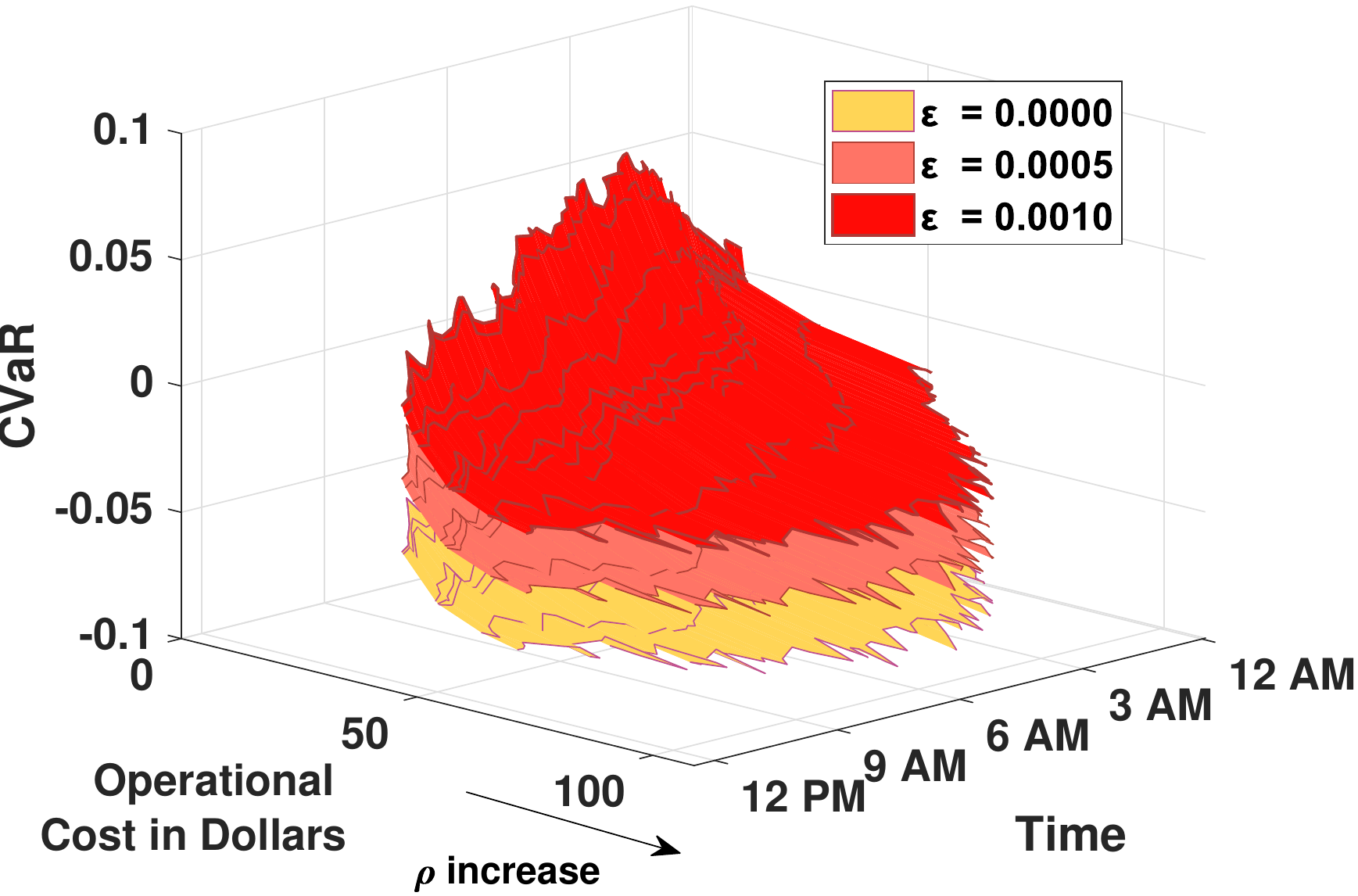}}
\subfigure{\label{fig:PVcurtail1} %% label for first subfigure
\includegraphics[width=1.7in]{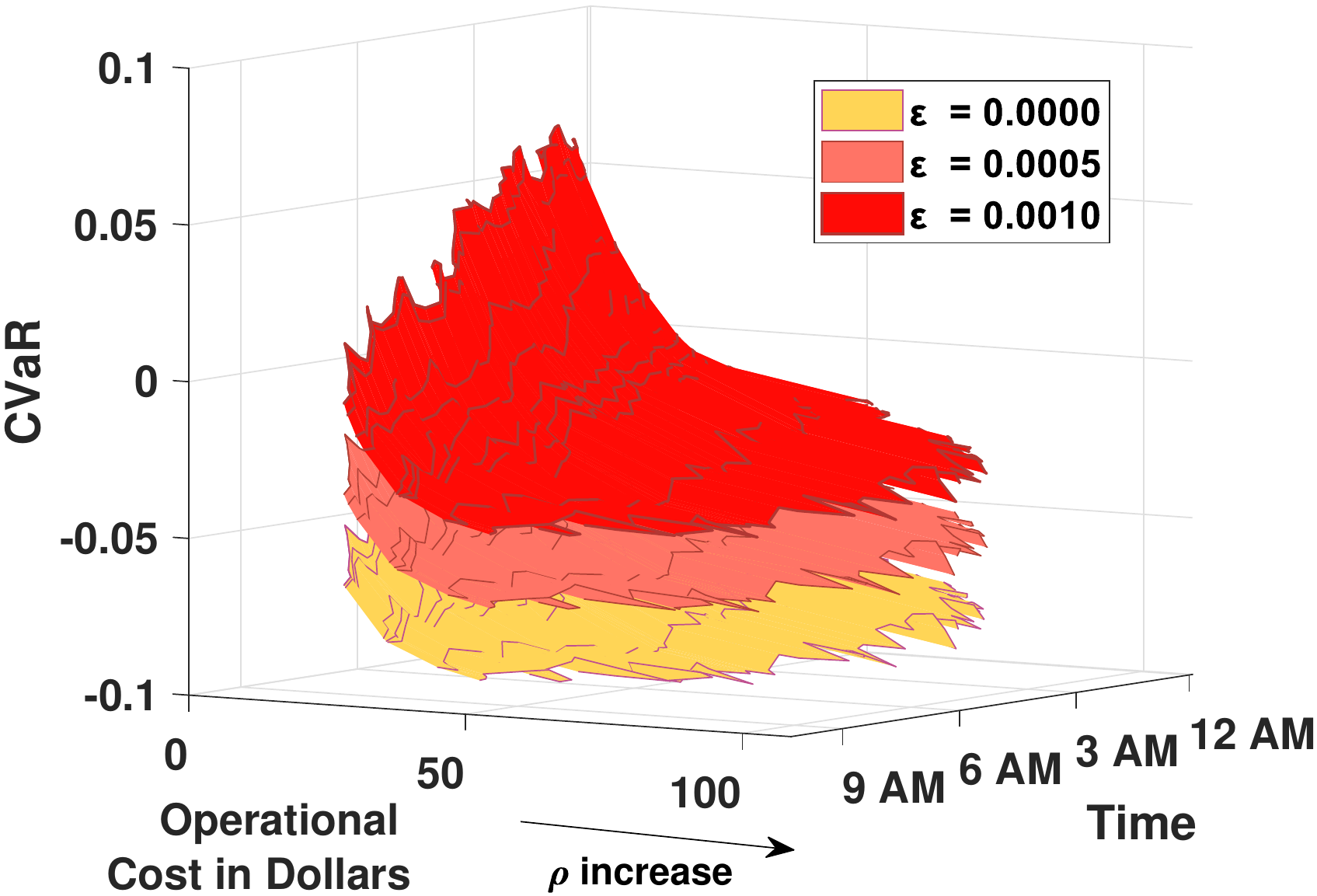}}

\caption{Tradeoffs between operational cost, conditional value at risk (CVaR) of  voltage constraints, and robustness to sampling errors for $\varepsilon = 0.0000, 0.0005, 0.0010$);  parameters $\rho$ and $\varepsilon$ are varied to tests  different weighting settings and radii of the Wasserstein ball, respectively. We present four views from different directions to avoid occlusion.}
\label{fidg:distributiontradeoff}
\end{figure*}

In our framework, there are two key parameters, $\rho$ and $\varepsilon$, that explicitly adjust trade offs between performance and constraint violation risk, and robustness to sampling errors. Fig. \ref{fidg:distributiontradeoff} illustrates the basic tradeoffs between operational cost and CVaR values of voltage constraint violations during a 24-hour operation for various values of $\rho$ and $\varepsilon$. It can be readily seen that as $\rho$ increases, operational cost increases, but CVaR decreases since the risk term is emphasized. Notice that with the increasing of $\varepsilon$, the \emph{estimated} risk is higher so that the solution is more conservative and leads to a lower risk of  constraint violation;  larger Wasserstein balls lead to higher robustness to sampling errors. 
These parameters offer system operators explicit data-based tuning knobs to systematically set the conservativeness of operating conditions.

\begin{figure*}[t!]
\centering
\subfigure[Total curtailed active power $\varepsilon = 0.0000$]{\label{fig:PVcurtail1} %% label for first subfigure
\includegraphics[width=2.2in]{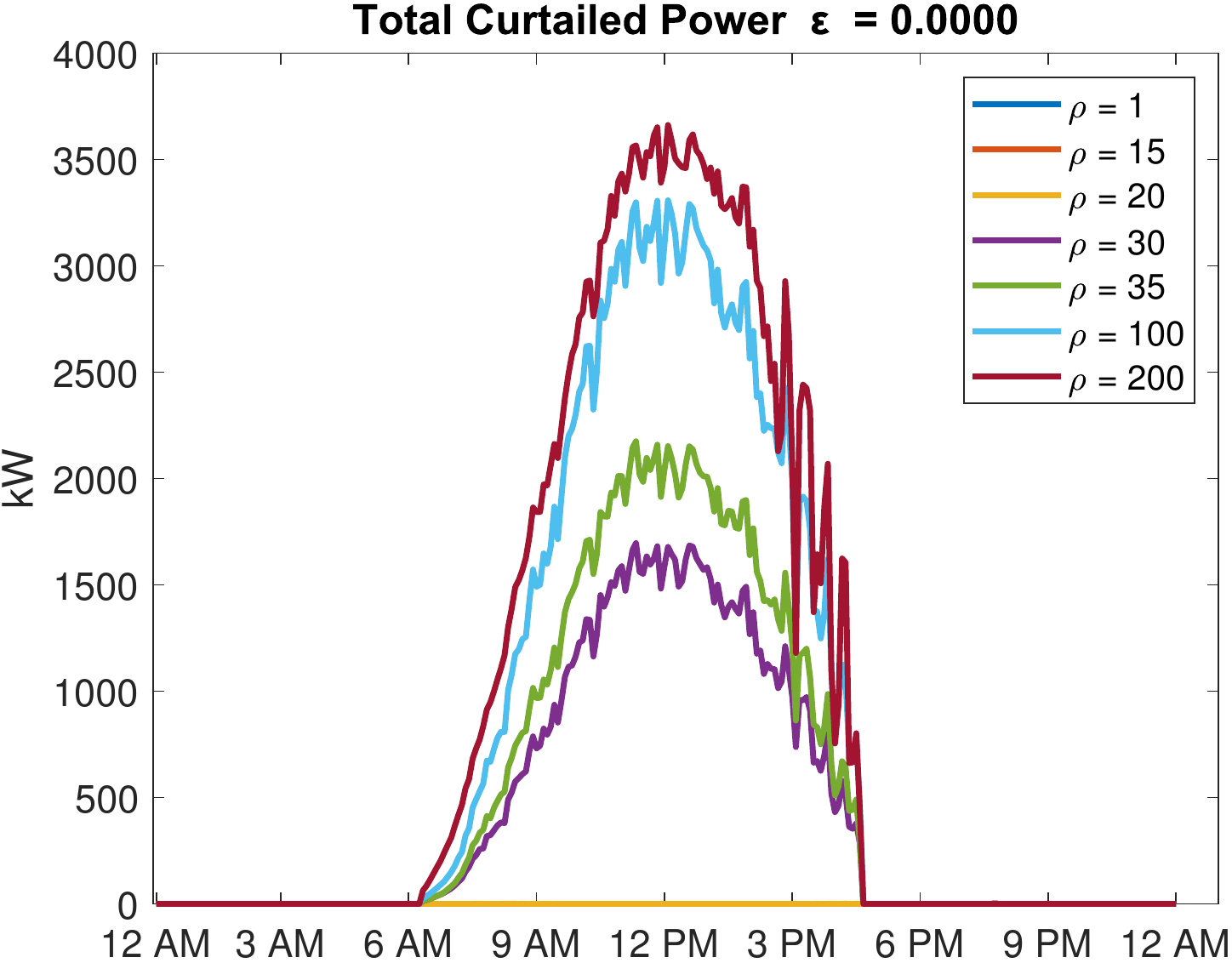}}
\subfigure[Total curtailed active power $\varepsilon = 0.0005$]{\label{fig:PVcurtial2} %% label for first subfigure
\includegraphics[width=2.2in]{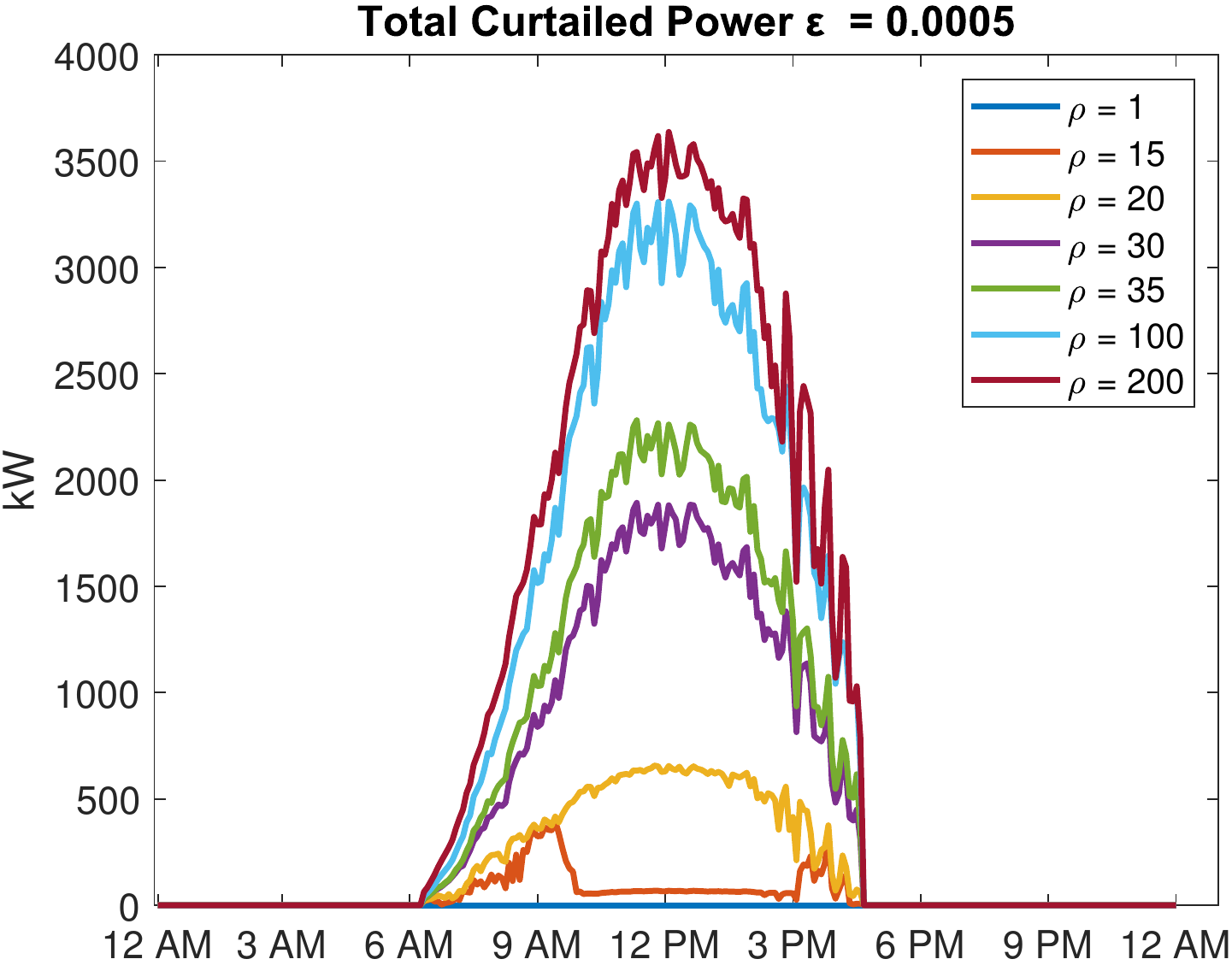}}
\subfigure[Total curtailed active power $\varepsilon = 0.0010$]{\label{fig:PVcurtail3} %% label for first subfigure
\includegraphics[width=2.2in]{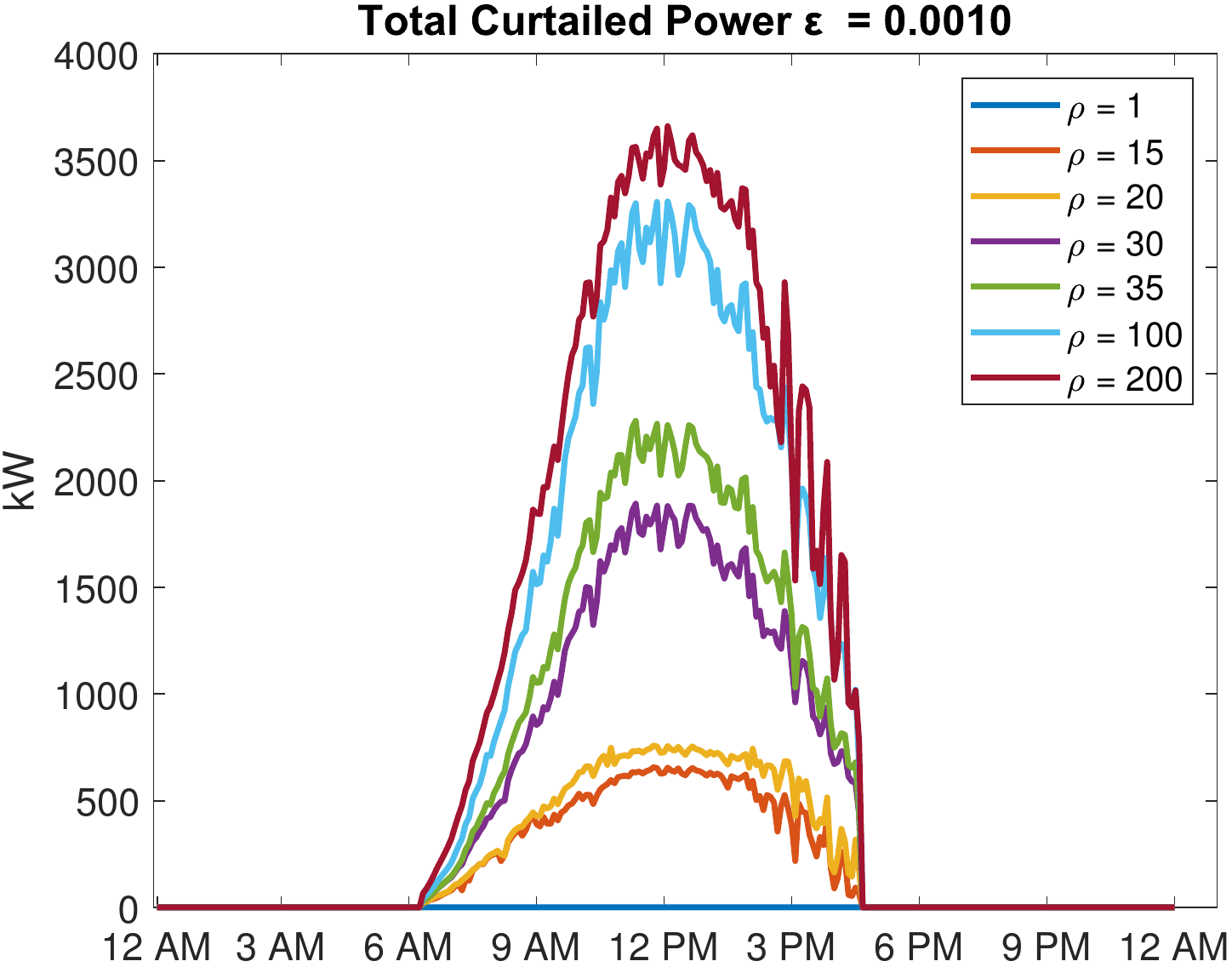}}
\subfigure[Power drawn from substation $\varepsilon = 0.0000$]{\label{fig:substation1} %% label for second subfigure
%%%%
\includegraphics[width=2.2in]{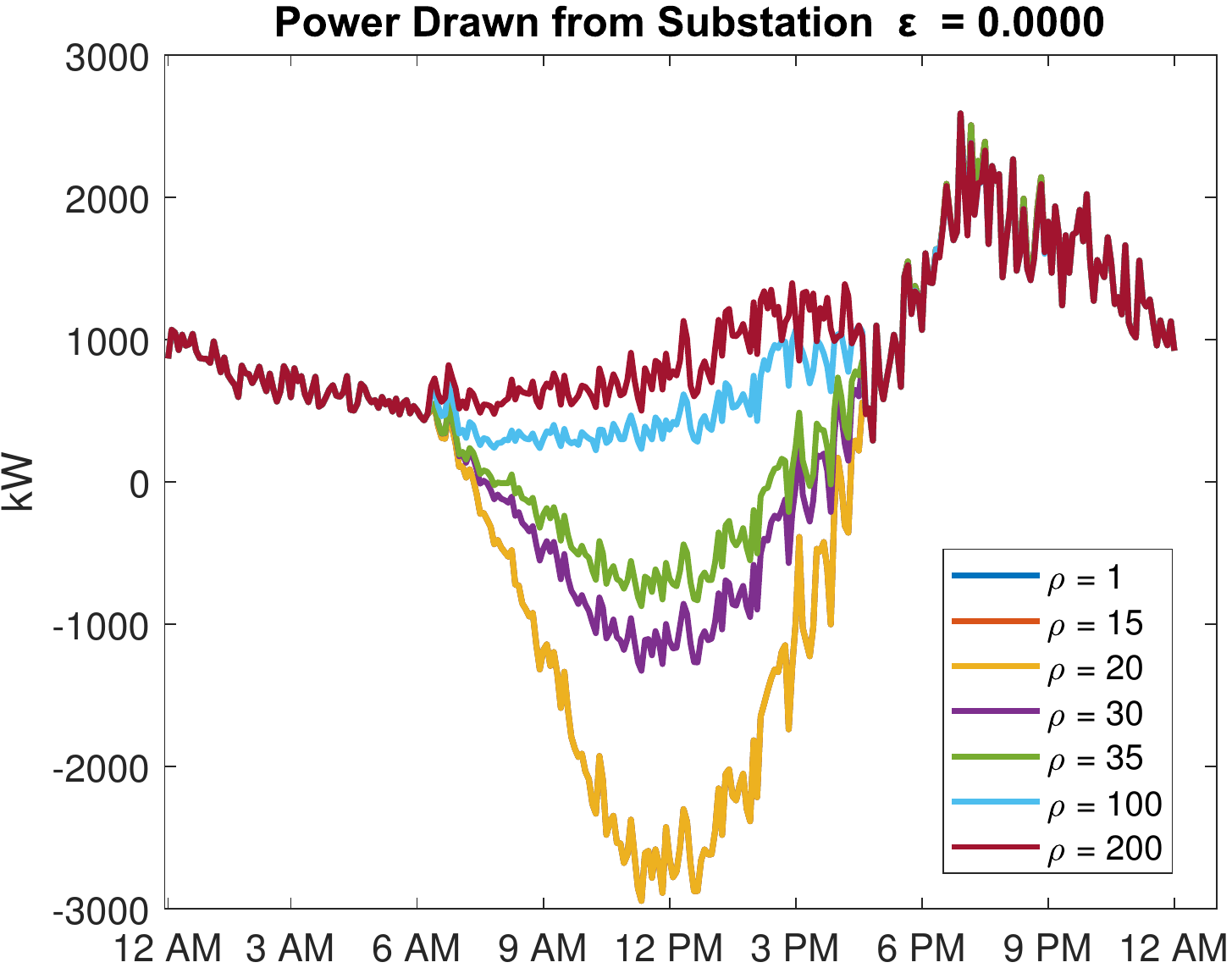}}
\subfigure[Power drawn from substation $\varepsilon = 0.0005$]{\label{fig:substation2} %% label for second subfigure
\includegraphics[width=2.2in]{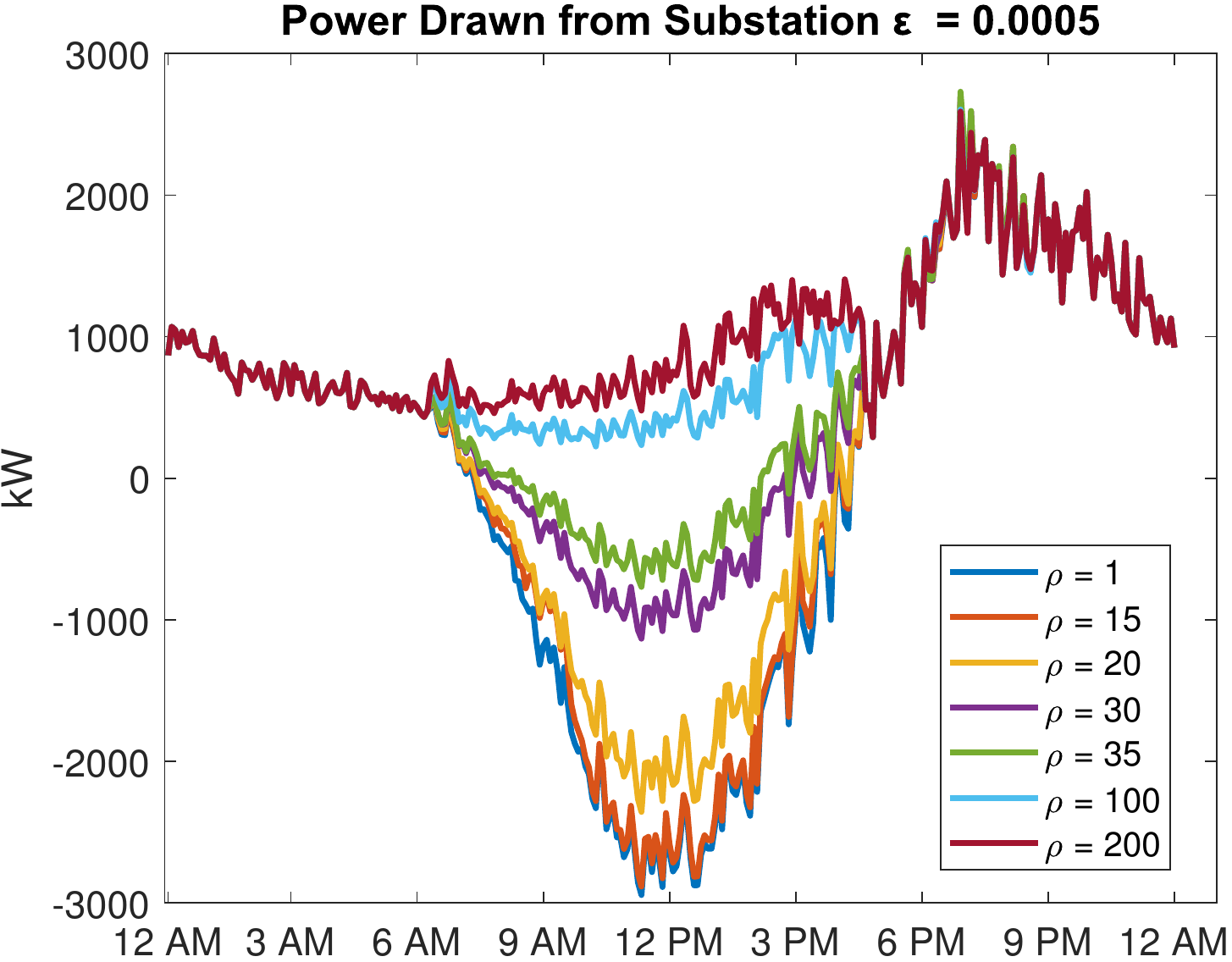}}
\subfigure[Power drawn from substation $\varepsilon = 0.0010$]{\label{fig:substation3} %% label for second subfigure
\includegraphics[width=2.2in]{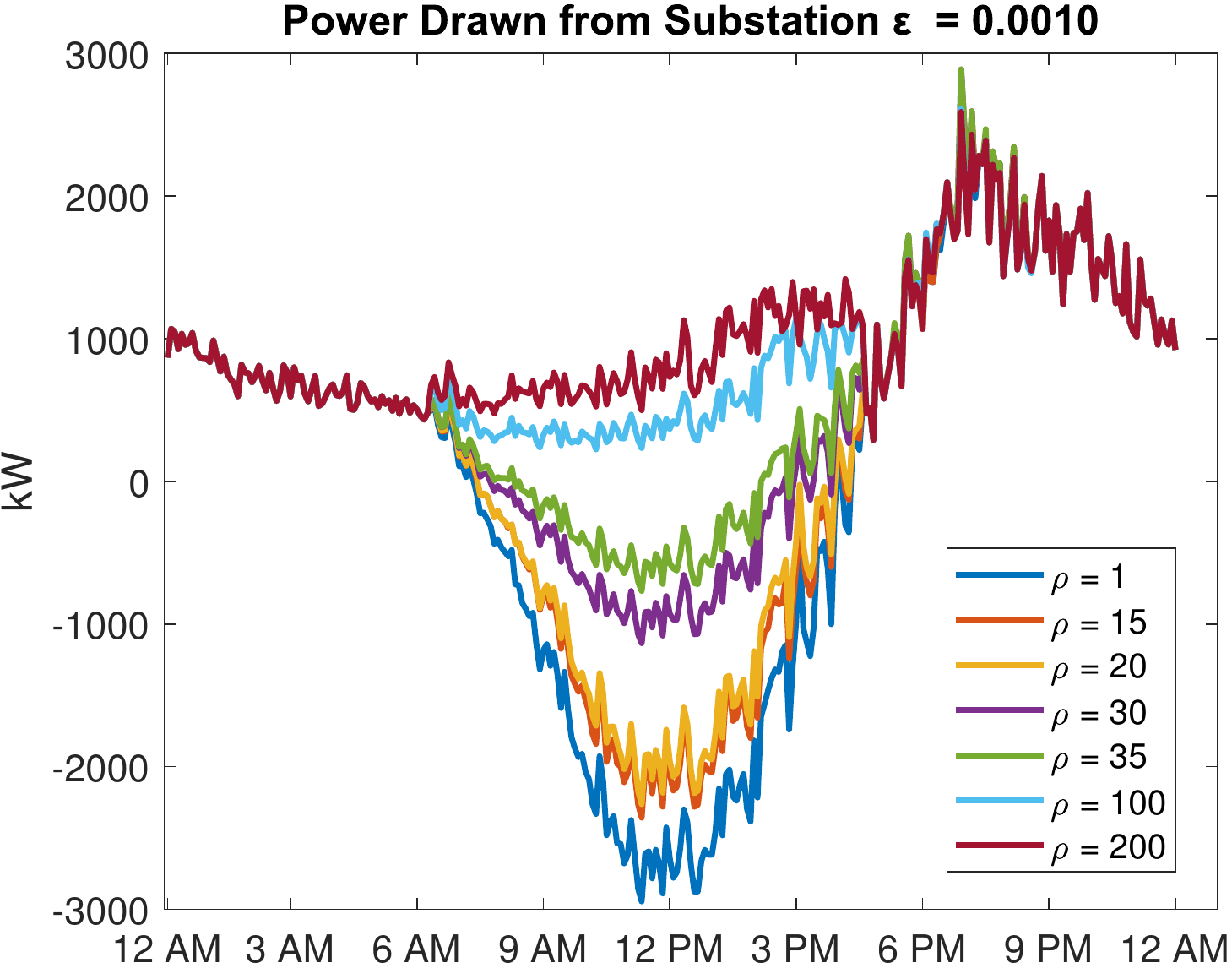}}

\caption{Comparison of active power curtailment and power purchased from substation for various values of risk aversion $\rho$ and Wasserstein radius $\varepsilon$. As these parameters increase, more active power from PV is curtailed and more power is drawn from the substation, leading to a lower risk of constraint violation and a higher operating cost.}
\label{fig:distributionPowerFlow} %% label for entire figure
\end{figure*}

Fig. \ref{fig:PVcurtail1}-\ref{fig:PVcurtail3} shows the aggregated solar energy curtailment and substation power purchases for varying risk aversion $\rho$ and Wasserstein radius $\varepsilon$. In order to prevent voltages over 1.05 p.u., the available solar energy must be increasingly curtailed as the risk aversion parameter $\rho$ increases. As a result, the network must import more power from the substation. The increasing curtailment of solar energy and purchase of power drawn from the substation lead to significantly higher operational cost. 

\begin{figure*}[hbtp!]
\centering
\includegraphics[width=7in]{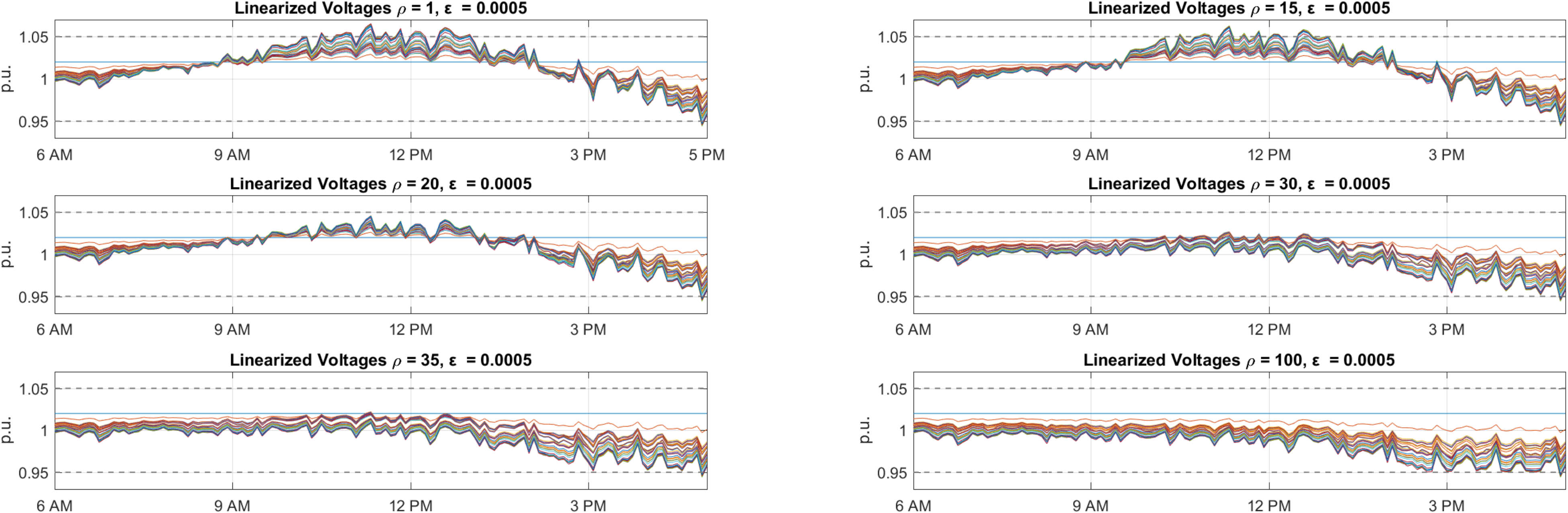}
\caption{Optimal network voltage profiles for varying $\rho$ and $\varepsilon$. Overvoltages are reduced as $\rho$ increases. }
\label{fidg:voltageprofile}
\end{figure*}

However, these decisions will also lead to more stable voltage profiles, as shown in Fig. \ref{fidg:voltageprofile}. When $\rho$ is small, there is almost no curtailment, causing overvoltages at several buses. As $\rho$ increases, more active power is curtailed, and all voltages move below their upper limit. Similar comments apply for varying the Wasserstein radius $\varepsilon$. For example, fixing $\rho$ and increasing $\varepsilon$ also results in more curtailment and lower voltage magnitude profiles, which leads to better robustness to solar energy forecast errors. 

\begin{figure*}[hbtp!]
\centering
\includegraphics[width=7in]{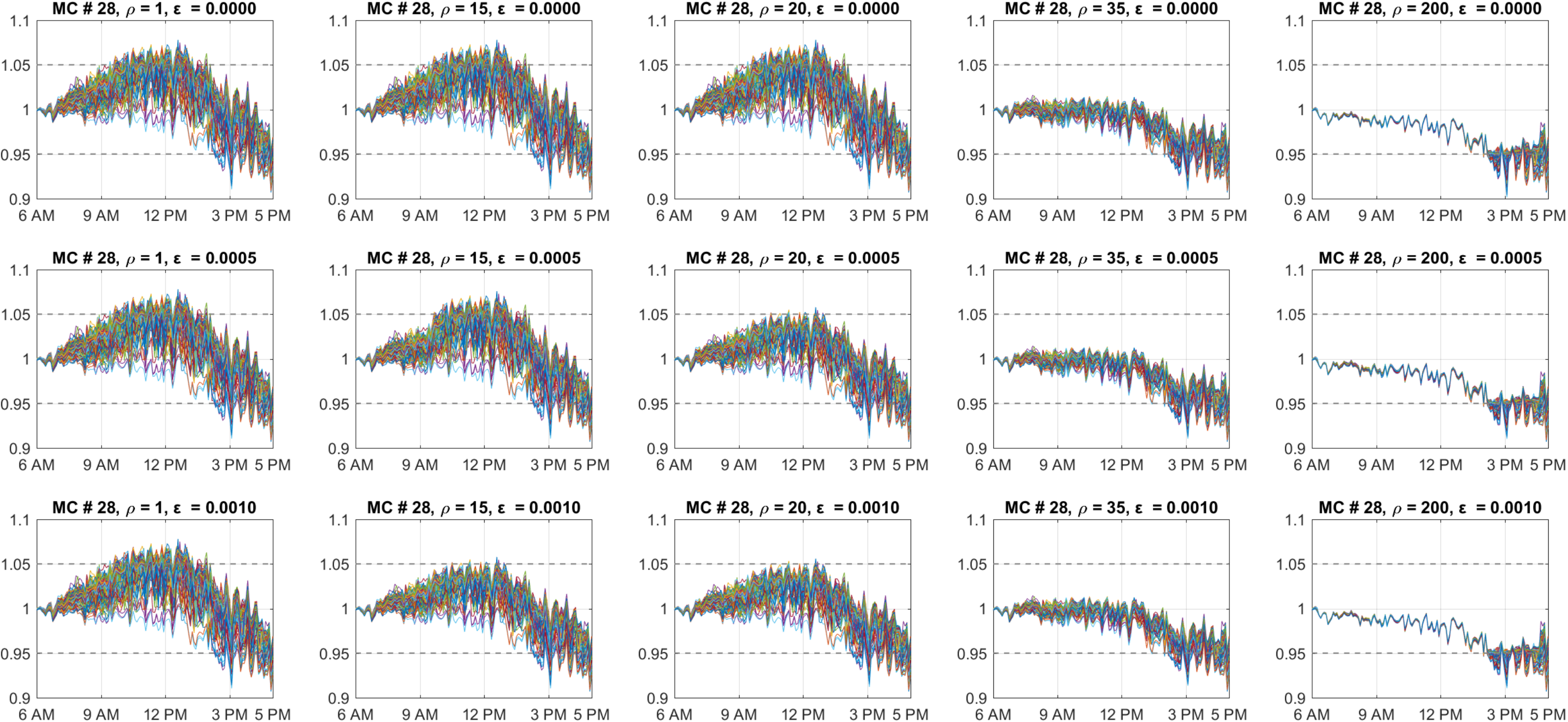}
\caption{Monte Carlo simulation results of the voltage profiles at node \#28 resulting from the full distributionally robust closed-loop model predictive control scheme. We validate that in closed-loop larger values of $\varepsilon$ and $\rho$ yield more conservative voltage profiles.  }
\label{fidg:MonteCarlosimulationDistribution}
\end{figure*}

Finally, we evaluate out-of-sample performance by implementing the full closed-loop distributionally robust MPC scheme over the 24 hour period with a 15 minute planning horizon. Monte Carlo simulations with 100 realizations of forecast errors over the entire horizon are shown in in Fig. \ref{fidg:MonteCarlosimulationDistribution}. We subsampled new solar energy forecast errors from the training dataset. The closed-loop voltage profiles based on MPC decisions for all scenarios at node $28$ are shown (other nodes with overvoltages show qualitatively similar results). Again, it is clearly seen that larger values of $\varepsilon$ and $\rho$ yield more conservative voltage profiles. 

In summary, we conclude that the proposed data-based distributionally robust stochastic OPF is able to systematically assess and control tradeoffs between the operational costs, risks, and sampling robustness in distribution networks. The benefits of the open-loop stochastic optimization problems are also observed in the closed-loop multi-period distributionally robust model predictive control scheme. 

% a graph for the distribution system test feeder.

% data-based stochastic opf in transmission systems.

\section{\emph{N}-1 security problem in transmission systems}
In this section, we apply the proposed methodology from Part I in a transmission system to handle \emph{N}-1 line flow security constraints. The basic DC power flow approximation and device modeling is discussed in Part I. Here, we also incorporate \emph{N}-1 security constraints and associated contingency reactions due to uncertain wind power injections.

\subsection{System model}
We consider a transmission system with $N_G$ generators (e.g., conventional thermal and wind) connected to bus subset $\mathcal{N}_G \subseteq \mathcal{N}$. There are $N_L$ loads, $N_l$ lines, and $N_b$ buses. The outages included for \emph{N}-1 security consist of tripping of any single lines, generators or loads, yielding $N_\textrm{out} = N_G + N_L + N_l$ possible outages. We collect the outages corresponding to a generator, a line and a load in different sets $\mathcal{I}_G$, $\mathcal{I}_l$, and $\mathcal{I}_L$. The outages in total are $\mathcal{I} = \{0\} \cup \mathcal{I}_G \cup \mathcal{I}_L \cup \mathcal{I}_l$, where $\{0\}$ indicates no outage. 

The formulation of a \emph{N}-1 security problem is based on the following assumptions: 1) the power flow equations are approximated with DC power flow, as described in Part I (Section V); 2) each wind farm is connected to a single bus of the network; 3) load forecasting is perfect; 4) a single line outage can cause multiple generator/load failures.

The objective of the data-based distributionally robust stochastic DC OPF is to determine an optimal reserve schedule for responding to the wind energy forecast errors while taking the network security constraints into account. We define $P_\textrm{mis}^j \in \mathbf{R}$ for all $j \in \mathcal{I}$ as the generation-load mismatch given by
\begin{equation}\nonumber
P_\textrm{mis}^j = \left\{\begin{array}{cl}
P_L^j - P_G^j, & \textrm{if}~~ j \in \{0\} \cup \mathcal{I}_l\\
+ P_L^j, & \textrm{if}~~ j \in \mathcal{I}_L\\
- P_G^j, & \textrm{if}~~ j \in \mathcal{I}_G\\
\end{array} \right. ,
\end{equation}
where $P_L^j, P_G^j \in \mathbf{R}$ denote the power disconnection corresponding to the outage $j \in \mathcal{I}$. Define $P_G \in \mathbf{R}^{N_G}$, and $P_L \in \mathbf{R}^{N_L}$ as nodal generation and load injection vectors. 
For the generator or load failures, the power disconnection $P_G^j$ or $P_L^j$ is corresponding to the components in the vector $P_G$ or $P_L$. In the case of line failures $j\in \mathcal{I}_l$, the power disconnection $P_G^j$ or $P_L^j$ is the sum of the power loss caused by multiple failures. If there is no power disconnection caused by a line outage, or if $j = 0$ (no outage happens), the power mismatch is set to zero: $P_{\text{mis}}^j = 0$.

To respond to contingencies, we can also define another reserve policy response matrix  $\mathbf{R}_{\textrm{mis},t}^{j,d} := [R_{1}^{j,d},\ldots,R_{t}^{j,d}]^\intercal \in \mathbf{R}^t$, so that the affine reserve policy becomes
\begin{equation} \label{reservepoliciesN1}
\begin{split}
\mathbf{u}_t^d = D_t^d \bm{\xi}_t + \mathbf{R}^{j,d}_{\textrm{mis},t} P_{\textrm{mis}}^j + e_t^d, \forall j\in \mathcal{I}, d = 1, \ldots, N_d.
\end{split}
\end{equation}
The general constraint risk function of the line flow in Part I (Section V, Equation (17c)) is then given by
% \begin{equation}
% \begin{split}
% \mathds{E}f_1 \bigg( \sum_{d=1}^{N_d} \widetilde{\Gamma}^i_d(r_d + G_d\xi_\tau + C_d(A_dx_{t0}^d + B_d\pi_\tau^d(\xi_\tau))) \\
% + R^{i,\tau}P_m^i - \bar{p} \bigg) \leq 0, ~\forall i \in\mathcal{I}, \tau \in \mathcal{T}_t, d = 1,...,N_d,
% \end{split}
% \end{equation}

\begin{equation}
\begin{split}
f\bigg( \sum_{d=1}^{N_d} \widetilde{\Gamma}^{t,j}_d \Big\{ r_d^t + G_d^t\bm{\xi}_t + C_d^t \big[ A^d_t x_{0}^d + B^d_t (D_t^d\bm{\xi}_t \\
+ \mathbf{R}^{j,d}_{\textrm{mis},t} P_{\textrm{mis}}^j + e_t^d ) \big] \Big\} - \bar{p}_t  \leq 0 \bigg), ~\forall j \in\mathcal{I},
\end{split}
\end{equation}
where $\widetilde{\Gamma}^{t,j}_d \in \mathbf{R}^{2Lt \times t}$ maps the power injection of each device in the case of $j$-th outage.

\subsection{Data-based stochastic OPF implementation}
We now use the modeling here and in Part I and the affine control strategy  \eqref{reservepoliciesN1} to formulate a data-based distributionally robust stochastic DC OPF for transmission systems that also incorporate \emph{N}-1 security constraints. The goal is to balance tradeoffs between cost of thermal generation, CVaR values of the line flow constraints, and sampling error robustness. The generation with reserve policy in the cost function is given by $P_G^{d,t} = \left[ D_t^d\xi_t + e_t^d\right]_t$. The operational cost of generators is
\begin{equation} \nonumber
J_\textrm{Cost}^t = \sum_{d \in \mathcal{N}_G}  c_{1,d} [P_G^{d,t}]^2 + c_{2,d} [P_G^{d,t}] + c_{3,d},
\end{equation}
which captures nominal and reserve costs of responding to wind energy forecast errors. The \emph{N}-1 security reserve cost is not included to simplify presentation, but this can also easily be included in our framework as an additional linear cost \cite{Vrakopoulou2}. 

With the proposed modeling in Part I (Section V), the updated data-based stochastic DC OPF is shown as follows. The decision variables are collected into $\mathbf{y}_t = \{\mathbf{D}_t, e_t, \mathbf{R}_{\textrm{mis},t}\}$. The random vector $\xi_t$ comprises all wind energy forecast errors.\\
\textbf{Data-based distributionally robust stochastic DC OPF}
\begin{subequations} \label{DBDCOPF}
\begin{align}
	& \nonumber \inf_{\mathbf{y}_\tau , \sigma_o^\tau} \sum_{\tau = t}^{t + \mathcal{H}_t} \bigg\{\mathds{E} [\hat{J}_{\textrm{Cost}}^\tau] +  \rho\sup_{\mathds{Q}_\tau \in \hat{\mathcal{P}}_\tau^{N_s}}\sum_{o=1}^{N_\ell} \mathds{E}^{\mathds{Q}_\tau}[\underline{\mathcal{Q}_o^\tau}]  \bigg\}, \\
	& = \inf_{\begin{subarray}{c} \mathbf{y}_\tau, \sigma_o^\tau, \\ \lambda_o^\tau, s_{io}^\tau, \varsigma_{iko}^\tau \end{subarray}} \sum_{\tau = t}^{\tau + \mathcal{H}_t} \bigg\{\mathds{E} [
\hat{J}_{\textrm{Cost}}^\tau] + \sum_{o=1}^{N_\ell} \bigg(\lambda_o^\tau\varepsilon_\tau + \frac{1}{N_s}\sum_{i=1}^{N_s} s_{io}^\tau\bigg) \bigg\},\\
	& \nonumber \textrm{subject to}\\
	& \label{TDRO1} \rho ( \underline{\mathbf{b}}_{ok}(\sigma_o^\tau) + \langle \underline{\mathbf{a}}_{ok}(\mathbf{y}_\tau), \hat{\xi}^{i}_\tau\rangle + \langle\varsigma_{iko}^\tau, \mathbf{d} -H\hat{\xi}^{i}_\tau\rangle ) \le s_{io}^\tau,\\ 
	& \label{TDRO2} \| H^\intercal\varsigma_{iko}^\tau-\rho \underline{\mathbf{a}}_{ok}(\mathbf{y}_\tau)\|_\infty \le \lambda_o^\tau,\\ 
	& \label{TDRO3} \varsigma_{iko}^\tau \geq 0,\\ 
	& \nonumber \frac{1}{N_s} \sum_{i=1}^{N_s}\sum_{d=1}^{N_d}\bigg[\widetilde{\Gamma}^{\overline{t},j}_d \Big\{r_d^{\overline{t}} + G_d^{\overline{t}}\bm{\hat{\xi}}_{\overline{t}}^i + C_d^{\overline{t}} \big[ A^d_{\overline{t}} x_{0}^d + B^d_{\overline{t}} (D_{\overline{t}}^d\bm{\hat{\xi}}_{\overline{t}}^i \\
	& \label{TN1Security}~~~~~~~~~~~~~~+ \mathbf{R}^{j,d}_{\textrm{mis},\overline{t}} P_{\textrm{mis}}^j + e_{\overline{t}}^d ) \big] \Big\} - \bar{p}_{\overline{t}} \bigg]_{[\underline{t},\overline{t}]} \leq 0, \\ 
	& \sum_{d=1}^{N_d}\Big[ (r_d^{\overline{t}} + C_d^{\overline{t}}(A^d_{\overline{t}}x_{0}^{d}+B^d_{\overline{t}}e_{\overline{t}}^d))\Big]_{[\underline{t},\overline{t}]} = 0,\\
	& \sum_{d=1}^{N_d}\Big[(G_d^{\overline{t}} + C_d^{\overline{t}}B^d_{\overline{t}}D^d_{\overline{t}})\Big]_{[\underline{t}, \overline{t}] } = 0, \\
	& \nonumber \forall i \le N_s, ~\forall j \in \mathcal{I}, ~\forall o \le N_\ell, ~k=1,2, ~\tau = t,\ldots, t+\mathcal{H}_t.
\end{align}
\end{subequations}
\subsection{Numerical results}
\begin{figure}[htbp!]
\centering
\includegraphics[width=3in]{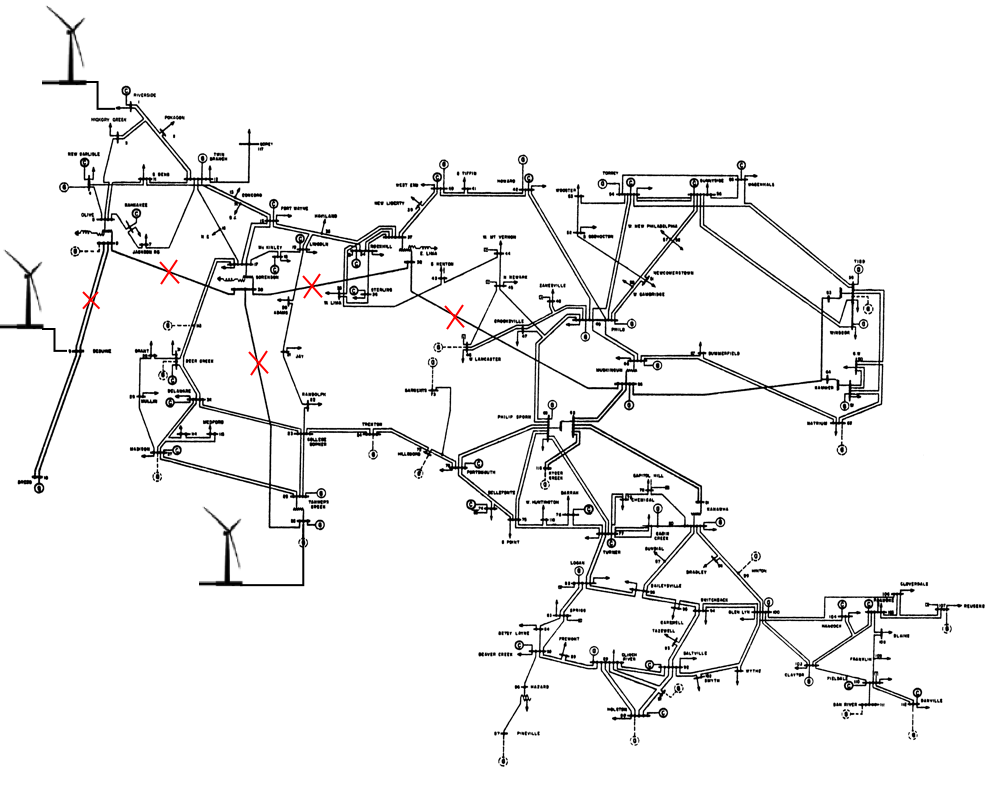}
\caption{IEEE 118-bus test network with multiple wind energy connections.}
\label{fidg:118transmission}
\end{figure}

We consider a modified IEEE 118-bus test system \cite{zimmerman2011matpower} to demonstrate our proposed data-based distributionally robust stochastic DC OPF shown in Fig. \ref{fidg:118transmission}. For simplicity, we only show results of a single-period stochastic optimization problem. As with the distribution network, it is straightforward to extend to multi-period closed-loop stochastic control using MPC.

\begin{figure*}[htbp!]
\centering
\subfigure[Predicted tradeoffs between operational cost and CVaR of line constraint violation.]{\label{fig:TransmissionsTradeoff} %% 
\includegraphics[width=3.37in]{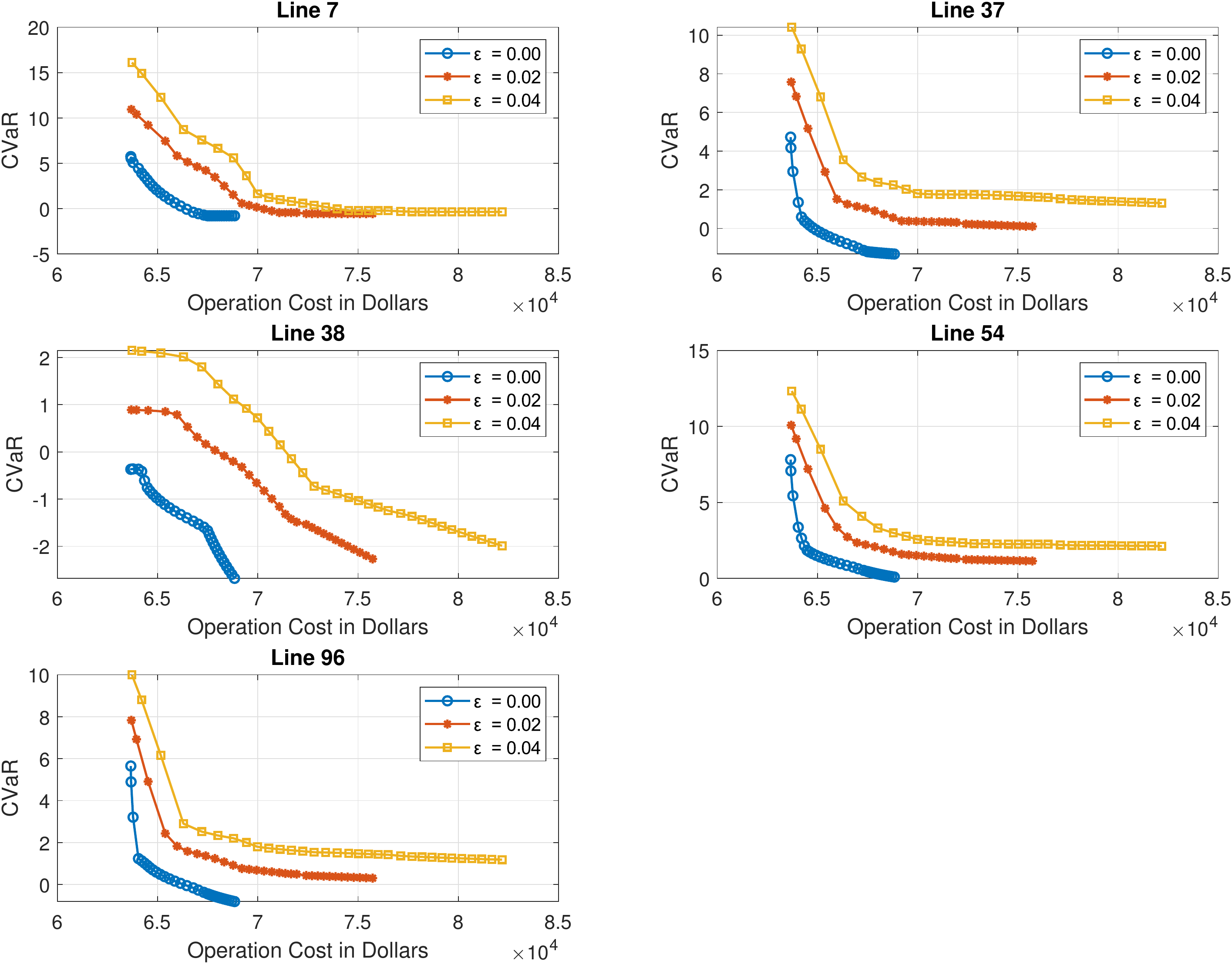}}
\subfigure[Out-of-sample performance is demonstrated by Monte Carlo simulation, with a controllable  level of  conservativeness. ]{\label{fig:TransmissionMonteCarloSimulations} 
\includegraphics[width=3.43in]{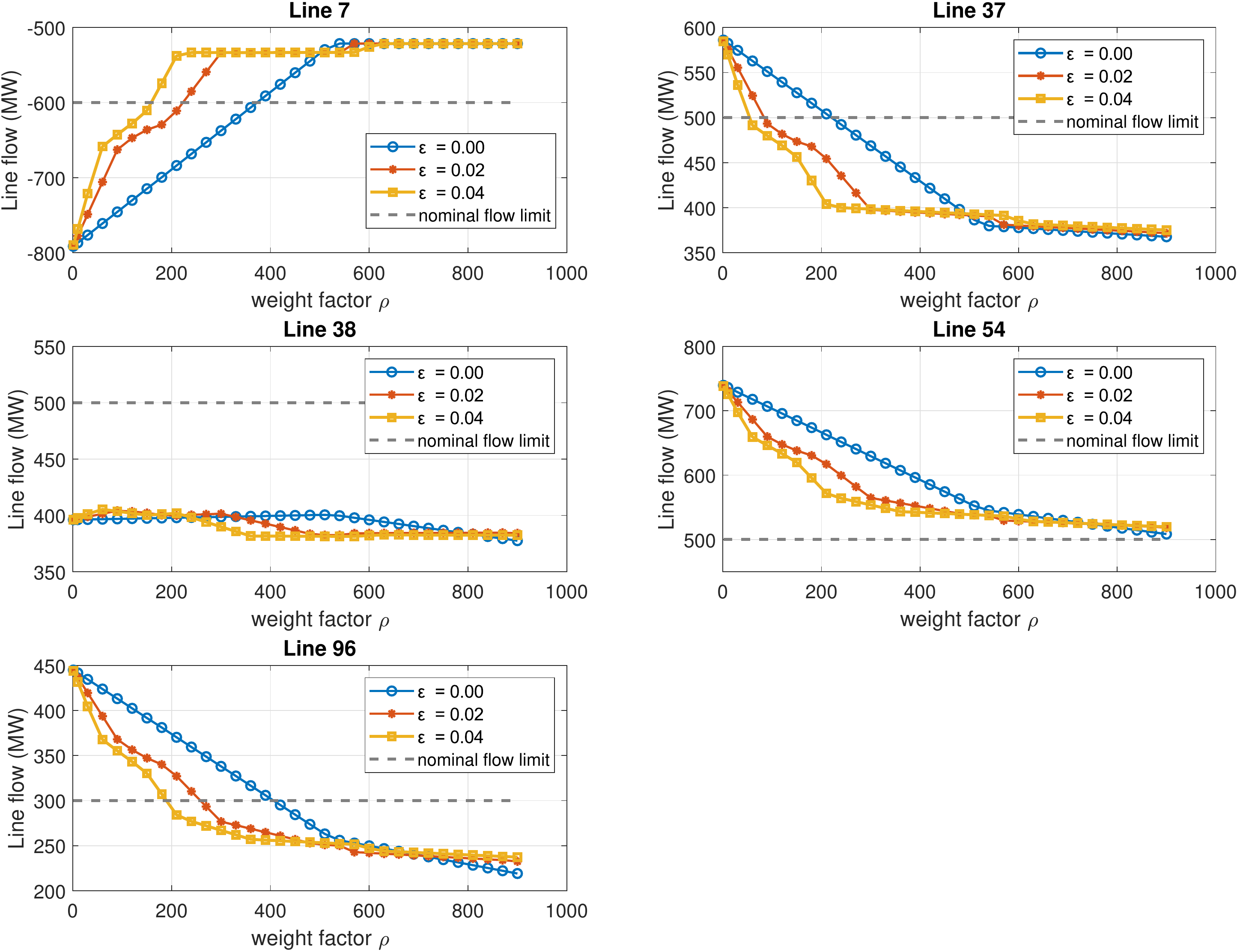}}
\caption{Five line flows are reported to illustrate the performance of the proposed data-based distributionally robust stochastic DC OPF in a transmission system.}
\label{fig:subfig} %% label for entire figure
\end{figure*}

\begin{figure*}[hbtp!]
\centering
\includegraphics[width=7in]{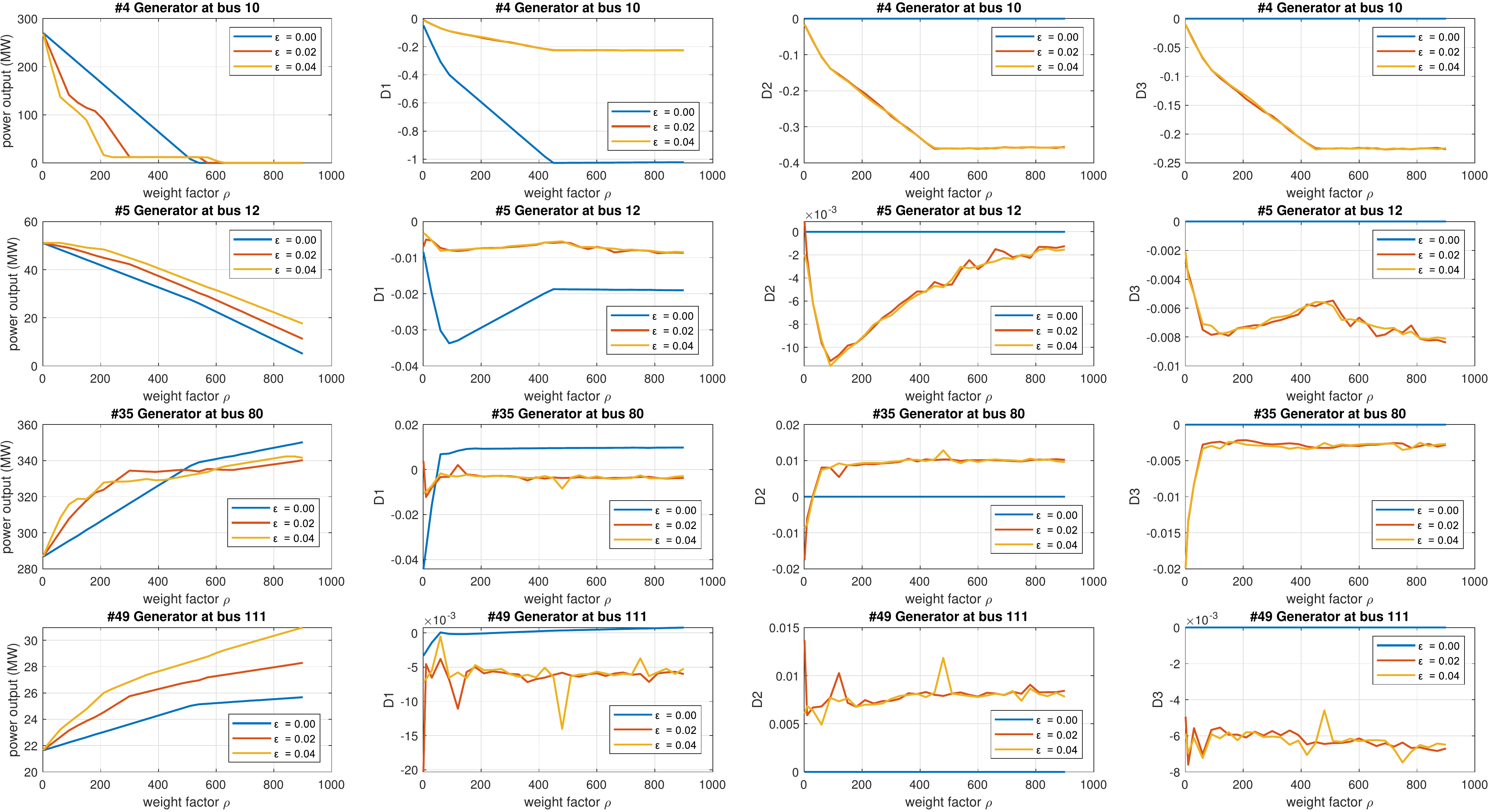}
\caption{Comparison of the coefficients of the  policies and output powers of selected generators for various values of risk aversion $\rho$ and Wasserstein metric $\varepsilon$. As these parameters increase, the risk of line flow constraint decreases at the expense of  higher operational costs.}
\label{fidg:ReservePolicesPowerOutputTransmission}
\end{figure*}

Three wind farms are connected to bus \#1, bus \#9, and bus \#26, with the normal feed-in power 500 MW, 500 MW, and 800 MW, respectively. The corresponding conventional generators at bus \#1 and bus \#26 are removed. The wind power forecast errors are derived from the real wind sampling data from hourly wind power measurements provided in 2012 Global Energy Forecasting competition (GEFCCom2012) \cite{wind}. The wind energy forecast errors are evaluated based on a simple persistence forecast, which predicts that the wind injection in the following period remains constant. It can be seen the forecast errors are highly leptokurtic, i.e, that the errors have significant outliers that make the distribution tails much heavier than Gaussian tails. We scale the forecasting errors to have zero mean and the standard deviations of 200 MW, 200 MW and 300 MW for the wind farms at bus \#1, \#9 and \#26, respectively. We consider five key lines, which deliver the wind power from the left side of the system to the right, as marked by the red crosses in Fig. \ref{fidg:118transmission}. The line flow limits are shown in Table \ref{table:LineParameters} and marked by gray dash lines in Fig.\ref{fig:TransmissionMonteCarloSimulations}.

\begin{table}[htbp!]
\centering
\caption{Five main channel lines data}
\begin{tabular}{|c|c|c|c|}
\hline
\# of line & From bus & To bus & Line flow limitation [MW]\\
\hline
7 & 8 & 9 & 600\\
\hline
37 & 8 & 30 & 500\\
\hline
38 & 26 & 30 & 500\\
\hline
54 & 30 & 38 & 500\\
\hline
96 & 38 & 65 & 300\\
\hline
\end{tabular}
\label{table:LineParameters}
\end{table}

% \begin{table}[htbp!]
% \centering
% \caption{Wind energy nominal deed-in power}
% \begin{tabular}{|c|c|c|}
% \hline
% \# bus & Nominal Feedin [MW] \\
% \hline
% 1 & 500\\
% \hline
% 9 & 500\\
% \hline
% 26 & 800\\
% \hline
% \end{tabular}
% \end{table}

To simplify our presentation, only these five lines flows are handled with distributionally robust optimization \eqref{TDRO1}-\eqref{TDRO3} in both directions; the remaining line flows are modeled by \emph{N}-1 security constraints \eqref{TN1Security} with nominal CVaR and sample average approximation (essentially equivalent to taking $\varepsilon = 0$), and no other local device constraint is included. Additionally, we assume no bounds on the wind power forecast errors $\xi_\tau$, hence $H$ and $\mathbf{d}$ in \eqref{TDRO1}-\eqref{TDRO2} are set to zero. It takes 10 seconds to solve \eqref{DBDCOPF} for each time step using the MOSEK solver \cite{mosek} via the MATLAB interface with CVX \cite{cvx} on a laptop with 16 GB of memory and a 2.8 GHz Intel Core i7.

Fig. \ref{fig:TransmissionsTradeoff} illustrates the solutions of the proposed data-based distributionally robust stochastic DC OPF for varying risk aversion and Wasserstein radius. Once again, the numerical results demonstrate fundamental tradeoffs between operational cost, CVaR values of line constraints violations, and robustness to sampling errors. The conservativeness of the generator policies are controlled by adjusting $\rho$ and $\varepsilon$. By explicitly using the forecast error training dataset and accounting for sampling errors, risks are systematically assessed and controlled. Since the forecast errors are non-Gaussian, existing methods may significantly underestimate risk \cite{guo2017stochastic}. Increasing $\varepsilon$ provides better robustness to sampling errors and guarantees out-of-sample performance. Note that CVaR values can be negative when the worst-case expected line flow is below the constraint boundary.

Fig. \ref{fig:TransmissionMonteCarloSimulations} illustrates the out-of-sample performance of the decisions via Monte Carlo simulation. For every  value of $\rho$ utilized to obtain the results of Fig. 8(a), we i) sampled (new) values from the training dataset, ii) implemented the decisions based on the solution of the problem, and iii) calculated the empirical line flows. The dashed gray line indicates the line flow limit. The number of scenarios for the  Monte Carlo simulations is 1000. From the simulation results, it is seen  that larger $\varepsilon$ ensures smaller line constraint violation for all lines except line 38. This happens because the risk objective is the sum of all five CVaRs, and a lower overall risk is achieved for certain values of $\rho$ and $\varepsilon$ by allowing higher risk of violating the flow limit of line 38. In general, it is possible to prioritize certain important constraints by weighting their associated risk higher compared to lower priority constraints. Again, Monte Carlo simulations demonstrate that conservativeness can be controlled explicitly by changing the Wasserstein radius $\varepsilon$ and the risk aversion parameter $\rho$. 

Fig.~\ref{fidg:ReservePolicesPowerOutputTransmission} illustrates the output powers and the coefficients of the reserve polices for selected generators for different values of the risk aversion $\rho$ and the Wasserstein radius $\varepsilon$. In order to satisfy the limit on the line flows, the scheduled power output of some generators (mostly located on the left side of the feeder or with cheaper cost profiles) are  reduced as the risk aversion parameter $\rho$ increase. As a result, some of the  generators (mostly located on the right side of the feeders or with high cost profiles) increase the power injection to supply the demand. With these settings, the risk of line flow constraint decreases as shown in Fig.~\ref{fig:TransmissionMonteCarloSimulations}, at the expense of  higher operational costs.  

% conclusion
\section{Conclusion}
We have illustrated the effectiveness and flexibility of our proposed multi-period data-based distributionally robust stochastic OPF methodology. We performed numerical experiments to balance overvoltages in distribution networks and $N-1$ security line flow risks in transmission networks. The flexibility of controllable devices was exploited to balance efficiency and risk due to high penetration renewable energy sources. In contrast to existing work, our method directly incorporates forecast error training datasets rather than making strong assumptions on the forecast error distribution, which allows us to leverage distributionally robust optimization techniques to achieve superior out-of-sample performance. Parameters in the optimization problems allow system operators to systematically select operating strategies that optimally trade off performance and risk. 

\section*{Supplementary materials}
Implementation codes for 1) data-based distributionally robust stochastic OPF and 2) data-based distributionally robust stochastic DC OPF can be download from \href{https://github.com/TSummersLab/Distributionally-robust-stochastic-OPF}{[Link]}. The general problem formulation and methodologies are presented in \cite{Guo2018DataDriven1}.

\section*{Acknowledgements}
The authors would like to thank Prof. Jie Zhang and Mr. Mucun Sun, at the Department of Mechanical Engineering, The University of Texas at Dallas, for their support on the solar forecasting.

\bibliography{reference2}
\bibliographystyle{ieeetr}
\end{document}